\documentclass[12pt]{article}

\usepackage{amsfonts,amsmath,amssymb,amsthm,graphicx,extarrows}
\usepackage[hyperfootnotes=false]{hyperref}
\usepackage{setspace} 

\def\proof{\noindent{\bf Proof:}\hskip10pt}        
\def\QED{\hfill $\Box$}

\font\tenmath=msbm10 scaled 1200
\font\sevenmath=msbm7 scaled 1200
\font\Fivemath=msbm5 scaled 1200
\newfam\mathfam \textfont\mathfam=\tenmath
\scriptfont\mathfam=\sevenmath \scriptscriptfont\mathfam=\Fivemath

\def \\ { \cr }
\def\R{\mathbb{R}}
\def \1{1 \mkern -6mu 1} 
\def\N{\mathbb{N}}
\def\E{\mathbb{E}}
\def\P{\mathbb{P}}

\def\R{\mathbb{R}}
\def \e{{\rm e}}
\def \d{{\rm d}}

\def \bC{{\bf C}}
\def \bY{{\bf Y}}
\def \bYabove{\overline{\bf Y}}
\def \bYbelow{\underline{\bf Y}}

\def \oS{\overline{S}}
\def \uS{\underline{S}}
\def \ob{\overline{\ensuremath{\beta}}}
\def \ub{\underline{\ensuremath{\beta}}}

\def \dt{{\ensuremath{\textup{d}}}}
\def \poi{{\ensuremath{\textup{Poi}}}}
\def \eqd{{\ensuremath{\stackrel{\textup{(d)}}{=}}}}

\newtheorem{theorem}{Theorem} 
\newtheorem{proposition}{Proposition} 
\newtheorem{lemma}{Lemma}
\newtheorem{corollary}{Corollary} 
\theoremstyle{definition}
\newtheorem{remark}{Remark}

\begin{document}

\title{Weak limits for the largest subpopulations 
  in Yule processes with high mutation probabilities}
\author{{Erich Baur\footnote{erich.baur@ens-lyon.fr} { and } Jean
    Bertoin\footnote{jean.bertoin@math.uzh.ch}}\\ ENS Lyon and
  Universit\"at Z\"urich}
\maketitle 
\thispagestyle{empty}

\begin{abstract}
  We consider a Yule process until the total population reaches size $n\gg
  1$, and assume that neutral mutations occur with high probability $1-p$
  (in the sense that each child is a new mutant with probability $1-p$,
  independently of the other children), where $p=p_n\ll 1$.  We establish a
  general strategy for obtaining Poisson limit laws for the number of
  subpopulations exceeding a given size and apply this to some mutation
  regimes of particular interest. Finally, we give an application to
  subcritical Bernoulli bond percolation on random recursive trees with
  percolation parameter $p_n$ tending to zero.
\end{abstract} 

\footnote{{\it Acknowledgment of support.} The research of the first author
  was supported by the Swiss National Science Foundation grant
  P300P2\_161011, and performed within the framework of the LABEX MILYON
  (ANR-10-LABX-0070) of Universit\'e de Lyon, within the program
  ``Investissements d'Avenir'' (ANR-11-IDEX-0007) operated by the French
  National Research Agency (ANR).}

\section{Introduction}
\label{Sintro}
We consider a system of branching processes with mutations specified as
follows. The underlying total population process is modeled by a standard
Yule process $Z$, that is a continuous-time birth process started from one
individual with unit birth rate per unit population size. We superpose
independent mutations, by declaring that a new-born child is a clone of its
parent with probability $p\in(0,1)$, and a mutant otherwise. Being a mutant
means that the individual obtains a new genetic type which was not present
before. We observe the process $Z$ at the instant when the $n$th individual
is born and group individuals of the same genetic type into subpopulations.

In this paper, we are interested in questions concerning the (asymptotic)
sizes of these subpopulations under strong mutations, the sense that
$p=p_n\rightarrow 0$. By approximating the population system from below and
above by two different processes, where sub-populations are independent and
have an explicit distribution, we develop a general strategy to obtain
non-trivial (Poisson) weak limits for the number of subpopulations
exceeding a given size (which might grow with $n$ as well).

We then discuss our strategy in the context of three qualitatively
different mutation regimes. For fixed $\ell\in\N$, we identify first
$p_n\sim an^{-1/\ell}$, $a>0$ fixed, as the regime in which, in the limit
$n\rightarrow\infty$, the largest subpopulations have size $\ell+1$. For
its number, we obtain a Poisson limit law and show that the number of
subpopulations of size $j$ for $j\in\{1,\ldots,\ell\}$ tends to infinity
(Theorem~\ref{thm:finitesubp} and Corollary~\ref{cor:1}).

Secondly, we discuss the regime $p_n\sim a/\ln n$. Since the size of the
subpopulation containing the ancestor is of order $n^{p_n}$, see
Proposition~\ref{prop:ancestralsubp}, this is the border-line case between
a bounded and an unbounded size for the ancestral subpopulation. We show
that the sizes of the largest subpopulations are concentrated around
$c_1\ln n +c_2\ln\ln n$, where $c_1$ and $c_2$ are positive constants
depending on $a$ (Theorem~\ref{thm:lognsubp2}). For the exact choice
$p_n=a/\ln n$ and $\lambda>0$, we find a correction
$c_3=c_3(a,\lambda)\in\R$ such that, with $y_n=c_1\ln n +c_2\ln\ln n +c_3$,
the number of subpopulations greater than $y_n$ converges along
subsequences $(y_{n(m)})\subset (y_n)$ with converging fractional part
to a Poisson$(\Lambda)$-distributed random variable, where $\Lambda$ is
expressed in terms of $a,\lambda$ and $(y_{n(m)})$ (Theorem~\ref{thm:lognsubp1}).

Thirdly, we study the case $1/\ln n\ll p_n\ll 1$. Here, it turns out that
the sizes of the largest subpopulations are to first order given by
$\e^{-1}p_n^{-1}n^{p_n}$. For $p_n\geq \ln\ln n/\ln n$ and given
$\lambda>0$, we compute a precise barrier  such that the number of
subpopulations exceeding this barrier follows in the limit the Poisson-law
with parameter $\lambda$ (Theorem~\ref{thm:greatsubp}).

This work originates from questions about Bernoulli bond percolation on
so-called random recursive trees, when their size $n$ tends to infinity and
the percolation parameter $p=p_n$ satisfies $p_n\rightarrow 0$. The
connection to branching processes stems from the fact that the genealogical
tree built from the first $n$ individuals in a standard Yule process can be
interpreted as a random recursive tree $T_n$ on $\{1,\ldots,n\}$.
Mutations in the Yule process can naturally be modeled on its genealogical
tree, by cutting the edges that connect mutants to its parent. Then the
connected subsets of vertices form the subpopulations of the same genetic
type. To put it differently, the connected components (clusters) on a
random recursive tree $T_n$ that arise from a Bernoulli bond percolation,
where each edge is erased with probability $1-p$ independently of each
other, can be viewed as the subpopulations in a Yule process with mutation
rate $1-p$, observed at the instant when there are $n$ individuals in total
in the system.
The strategy we develop here in terms of Yule processes allows a concise
analysis of cluster sizes, for any choice of $p_n$ tending to zero. For
sequences of $p_n$ such that $p_n\rightarrow 1$ or $p_n=p\in(0,1)$ remains
constant, similar connections between systems of branching processes and
percolation on increasing tree families have been utilized before in,
e.g.,~\cite{Ber,BeBr,Ba,BaBe}. The precise definition of a random recursive
tree, its connection to Yule processes and more references to existing
results on percolation will be discussed in
Section~\ref{sec:app-percolation}.

The rest of this paper is organized as follows. In the next
Section~\ref{sec:Yule}, we properly define the population system and
provide some heuristics for regimes of interest and Poisson limits. Then,
in Section~\ref{sec:strategy}, we explain our strategy for obtaining
Poisson limit laws for the number of subpopulations (or clusters) greater
than a given size. Section~\ref{sec:limitresults} contains our main
results; we exemplify our strategy by proving limit results for certain
mutation rates of particular interest. In the last
Section~\ref{sec:app-percolation}, we establish the link to percolation on
random recursive trees. Appendix~\ref{sec:appendix} contains some
(standard) estimates on Yule processes, which we use in our analysis.

{\bf Notation:} We let $\N=\{1,2,\ldots\}$. If $(a_n:n\in\N)$,
$(b_n:n\in\N)$ are two sequences of real numbers, we write $a_n\sim b_n$ if
$a_n/b_n\rightarrow 1$ as $n\rightarrow\infty$, and we write $a_n \ll b_n$
or $b_n \gg a_n$ if and only if $a_n/b_n \to 0$ as $n\to\infty$.  Moreover,
if $f(n)$ and $g(n)$ are two positive functions, we say $g(n)=O(f(n))$ if
there exists $M>0$ such that $g(n)\leq Mf(n)$ for all $n\in\N$, and we
write $f(n)=o(g(n))$ if $g(n)/f(n) \rightarrow 0$ as $
n\rightarrow\infty$. We will use the letters $c$ or $C$ for small or large
generic constants that do not depend on $n$. Their values may change from
line to line.
\section{Yule processes with mutations}
\label{sec:Yule}
In this section, we introduce the population system we work with, and
present some basic results and heuristics. For background on Yule
processes, we refer to Appendix~\ref{sec:appendix}, and for more
information to Chapter $3$ of~\cite{AtNe}.

\subsection{A population system with infinitely many types}
Let $Z=(Z(t):t\geq 0)$ be a standard Yule process with unit birth rate per
individual and started from one single particle, i.e. $Z(0)=1$. Let
$p\in[0,1]$. We superpose mutations on $Z$ as follows: When a new child is
born, we declare it to be a clone of its parent with probability $p$, and a
mutant with a new genetic type different from all the previous types with
probability $1-p$.  We let $b_1^{(p)}=0$ and write
$0<b_2^{(p)}<b_3^{(p)}<\dots$ for the sequence of consecutive birth times
of individuals which are mutants. More specifically, the ancestor present
at time $b_1^{(p)}=0$ has genetic type $1$, as well as its clone
children. The first mutant is born at time $b_2^{(p)}$ and has genetic type
$2$. The next mutant appears in the system at time $b_3^{(p)}$ and receives
genetic type $3$; it is a mutant of an individual of either type $1$ or of
type $2$, and so on. We group individuals of the same genetic type into
subpopulations, so that at time $t$, the population system consists of at
most $Z(t)$ many subpopulations. In contrast to~\cite{BeBr}, we will here
not be interested in the genealogical structure, but only in the sizes of
the subpopulations.

For $i\in\N$, we let $(Y^{(p)}_i(t):t\geq 0)$ denote the subpopulation
process counting the individuals of type $i$, with $Y^{(p)}_i(t)=0$ if
$t<b_i^{(p)}$.  It should be clear that the processes
$(Y^{(p)}_i(b_i^{(p)}+t):t\geq 0)$, $i\in\N$, form independent Yule
processes with birth rate $p$ per individual and started from one single
individual each. The number of subpopulations of different genetic types
present at time $t\geq 0$ is denoted by
$$
T^{(p)}(t)=\max\left\{i\in\N_0:b_i^{(p)}\leq t\right\}.
$$
Viewed as a process in $t$, $T^{(p)}=(T^{(p)}(t),t\geq 0)$ is a counting
process started from $1$ which grows at rate $(1-p)Z(t)$. Its predictable
compensator is absolutely continuous with density $(1-p)Z(t)$,
that is
$$
T^{(p)}(t)-\int_0^t(1-p)Z(s)\dt s,\quad t\geq 0,
$$
is a martingale. See, e.g., \cite[Theorem 9.15]{Kl}. 

We next build a process $\bY^{(p)}=(\bY^{(p)}(t):t\geq 0)$ by setting
$$
\bY^{(p)}(t)=\left(Y_1^{(p)}(t),Y_2^{(p)}(t),\ldots\right),\quad t\geq 0.
$$
Note that $\bY^{(p)}(0)=(1,0,\ldots)$. Clearly, we can retrieve the total
population size $Z(t)$ at time $t$ from $\bY^{(p)}$,
$$
Z(t)=\sum_{i=1}^\infty Y_i^{(p)}(t).
$$
It follows from its construction that the process $\bY^{(p)}$ is Markovian
with transition rates at time $t\geq 0$ for $y_1,\ldots,y_k\in\N$ given by
\begin{align*}
  (y_1,\ldots,y_k,0,\ldots) &\mapsto (y_1,\ldots,y_k,1,0,\ldots)\quad \textup{at rate }(1-p)(y_1+\ldots+y_k),\\
  (y_1,\ldots,y_k,0,\ldots)&\mapsto(y_1,\ldots,y_{i-1},y_i+1,y_{i+1},\ldots,y_k,0,\ldots)
  \quad \textup{at rate }py_i.
\end{align*}
Mostly we will stop $\bY^{(p)}$ at the instant when the $n$th particle is born,
$$
\tau_n=\inf\left\{t\geq 0: Z(t)=n\right\}.
$$
Obviously, we have $\tau_i\leq b_i^{(p)}$ for all $i\in\N$ and all $p\in[0,1]$.  Moreover, since
$\e^{-t}Z(t)$ is a martingale which converges almost surely to a standard
exponential $\mathcal{E}$ as $t$ tends to infinity, see Appendix~\ref{sec:appendix},
\begin{equation}
\label{eq:behaviortau_n}
\lim_{n\rightarrow\infty}\left(\tau_n-\ln n\right) =-\ln
\mathcal{E}\quad\hbox{almost surely}.
\end{equation} 

\subsection{The size of the ancestral subpopulation}
In this section, we point at a simple limit theorem in distribution for 
the size $Y_1^{(p)}(\tau_n)$ of the subpopulation at time $\tau_n$ having the
same genetic type as the ancestor individual. We obtain the following
characterization when $p=p_n\ll 1$. 

\begin{proposition}
\label{prop:ancestralsubp}
  For $\lambda >0$, denote by \textup{Geo}$(\lambda)$ a geometrically distributed
  random variable of parameter $\lambda$, and by \textup{Exp}$(1)$ a standard
  exponential random variable. Then the following holds:
\begin{enumerate}
\item If $0\leq p_n\ll 1/\ln n$, then $Y_1^{(p_n)}(\tau_n)\xrightarrow[n
  \to \infty]{(d)} 1.$
\item If $p_n\sim a/\ln n$ for some $a>0$, then $Y_1^{(p_n)}(\tau_n)\xrightarrow[n
  \to \infty]{(d)} \textup{Geo}(\e^{-a}).$
\item If $1/\ln n\ll p_n\ll 1$, then $n^{-p_n}Y_1^{(p_n)}(\tau_n)\xrightarrow[n
  \to \infty]{(d)} \textup{Exp}(1).$
\end{enumerate}
\end{proposition}
\proof Assume $p_n\rightarrow 0$, and let $G_n=\left\{|\tau_n-\ln n|\leq
p_n^{-1/2}\right\}$. From~\eqref{eq:behaviortau_n}, we know $\P(G_n)\rightarrow 1$ as
$n\rightarrow\infty$. For $(a)$ and $(b)$, it thus suffices to
show that for $x\in\N$, 
$$
\P\left(Y_1^{(p_n)}(\tau_n)>x;\, G_n\right)\rightarrow
1-F(x)\quad\textup{as }n\rightarrow\infty,
$$
where $F$ is the distribution function of the stated limit in $(a)$ or
$(b)$, respectively. Since
$Y_1^{(p_n)}(t)$ is monotone increasing in $t$, we have 
\begin{equation}
\label{eq:ancestral-eq}
\P\left(Y_1^{(p_n)}(\ln n -p_n^{-1/2})>x;\, G_n\right)\leq
\P\left(Y_1^{(p_n)}(\tau_n)>x;\, G_n\right)\leq \P\left(Y_1^{(p_n)}(\ln n
  +p_n^{-1/2})>x\right).
\end{equation}
Since $Y^{(p_n)}(t)$ has a geometric law with success probability
$\exp(-p_nt)$, we obtain 
$$
\P\left(Y_1^{(p_n)}(\ln n
  +p_n^{-1/2})>x\right)=\left(1-\exp\left(-\left(p_n\ln n+p_n^{1/2}\right)\right)\right)^{x}.
$$
It is readily checked that when $x\in\N$ and $n\rightarrow\infty$, the right side
converges for $p_n\ll 1/\ln n$ to zero, whereas for $p_n\sim a/\ln n$, it
converges to $(1-\e^{-a})^x$. Writing 
$$\P\left(Y_1^{(p_n)}(\ln n -p_n^{-1/2})>x;\, G_n\right) = \P\left(Y_1^{(p_n)}(\ln n -p_n^{-1/2})>x\right)+\P(G_n^c),$$
we see that the same convergence holds for the probability on the left side
in~\eqref{eq:ancestral-eq}. This shows $(a)$ and $(b)$. The proof of $(c)$
is entirely similar and left to the reader.  \QED
\begin{remark}
\begin{enumerate}
\item The above statements remain true $Y_1^{(p)}(\tau_n)$ is replaced by any
  $Y_i^{(p)}(\tau_n)$ for $i\in\N$ fixed.
\item For $p\in (0,1)$ fixed, it has been shown independently in~\cite{BaBe} and \cite{Mo} that
  $n^{-p}Y_1^{(p)}(\tau_n)$ converges to a Mittag-Leffler
  distributed random variable with parameter $p$. In turn, the latter converges to a
  standard exponential random variable when $p\rightarrow 0$. If
  $p_n\rightarrow 1$ as $n\rightarrow\infty$, then
  $n^{-p_n}Y_1^{(p_n)}(\tau_n)\rightarrow 1$ a.s., see~\cite{Ba}.
\item We note that for each $p\in [0,1]$, the mean of $Y_1^{(p)}(\tau_n)$ can be
computed exactly. Indeed, we have
$$
\E\left[Y_1^{(p)}(\tau_n)\right] = \sum_{i=1}^n\P\left(i\in
    Y_1^{(p)}(\tau_n)\right) = \sum_{i=1}^n\E\left[p^{D(i)}\right],
$$
where $D(i)$ is the insertion depth of $i$, i.e., $D(i)$ is the height of $i$
in the genealogical tree of the Yule process stopped at time $\tau_n$ (with
$D(1)=0$). From~\cite{DoSm} we know
that $D(i) \eqd \sum_{j=1}^{i-1}\xi_j$, where $\xi_j$ are independent
Bernoulli random variables of parameter $1/j$. Consequently,
\begin{align*}
\sum_{i=1}^n\E\left[p^{D(i)}\right] =
1+p\sum_{i=2}^n\prod_{j=2}^{i-1}\E\left[p^{\xi(j)}\right]
&=1+p\sum_{i=2}^n\prod_{j=2}^{i-1}\left(\frac{p}{j}+ 1-\frac{1}{j}\right)\\
&=1+p\sum_{i=1}^{n-1}\frac{1}{i!}\prod_{j=2}^{i}(j-1+p).
\end{align*}
\end{enumerate}
\end{remark}

\subsection{Poisson heuristic for the number of large subpopulations}
\label{sec:poisson-heu}
Clearly, when $p_n\ll 1/n$, the $n$ first individuals are all mutants and
all subpopulations at time $\tau_n$ are singletons with high
probability. When $p_n\sim a/n$, standard Poisson approximation to the
binomial law (and the fact that the genealogical tree of $Z(\tau_n)$ is not
star-like) shows that the number of subpopulations of size $2$ is
asymptotically Poisson$(a)$-distributed. For general $\ell > 1$, it is
however {\it a priori} not obvious in which regime subpopulations of size
$\ell+1$ will emerge, and if so, how many of them. Following the tradition
of Aldous \cite{Aldous}, we shall now present an heuristic argument to
determine this regime, and will check later on that our informal approach
can actually be made rigorous.

We are interested in the number of subpopulations of size
$>\ell $ at time $\tau_n$,
$$
N^{n,p}(\ell)=\sum_{i=1}^{\infty}\1_{\left\{Y^{(p)}_i(\tau_n)>\ell\right\}}.
$$

For $i\gg 1$ and $p\ll 1$, the birth time $b^{(p)}_i$ of the $i$th mutant
is close to the birth time $\tau_i$ of the $i$th individual, and by
\eqref{eq:behaviortau_n}, we have $b^{(p)}_i-\ln i \sim -\ln \mathcal{E}$
a.s. Thus the distribution of $Y^{(p)}_i(\tau_n)$ should be close to that
of a Yule process with birth rate $p$ evaluated at time
$\tau_n-b^{(p)}_i\sim \ln(n/i)$, that is to a geometric distribution with
parameter $\exp(-p\ln(n/i))=(i/n)^p$. Therefore, the Bernoulli variable
$\1_{\left\{Y^{(p)}_i(\tau_n)>\ell\right\}}$ has parameter
$$\P(Y^{(p)}_i(\tau_n)>\ell)\sim \left(1-(i/n)^p\right)^\ell \sim \left( p \ln\frac{n}{i}\right)^\ell.$$

Let us now recall Le Cam's inequality for the reader's convenience.

\noindent{\bf Inequality of Le Cam~\cite{LC}} {\it Let $(a_n)_n\subset \N$
  be a sequence of integers. For each $n\in\N$, let $\xi_{n,i}$, $1\leq
  i\leq a_n,$ be independent Bernoulli random variables with parameters
  $q_{n,i}$. Set $S_n = \xi_{n,1} + \cdots + \xi_{n,a_n}$, and let
  $\lambda_n = \E[S_n]=\sum_{i=1}^{a_n} q_{n,i}$. Then
$$
\sum_{j=1}^\infty\left|\P\left(S_n=j\right)-\frac{\lambda_n^j\e^{-\lambda_n}}{j!}\right|\leq
  2\sum_{i=1}^{a_n} q_{n,i}^2.
$$
In particular, if $\lim_{n\rightarrow\infty}\lambda_n=\lambda\in(0,\infty)$
and $\lim_{n\rightarrow\infty}\sum_{i=1}^{a_n} q_{n,i}^2=0$, then
$S_n\xrightarrow[]{(d)} \poi(\lambda)$, where $\poi(\lambda)$ is a
Poisson$(\lambda)$-distributed random variable.}

Even though the variables $Y^{(p)}_i(\tau_n)$, $i\in\N$, which describe the
sizes of the subpopulations, are clearly not independent (obviously, they
add up to $n$), let us pretend for the purpose of this section that
$\1_{\left\{Y^{(p)}_i(\tau_n)>\ell\right\}}$ for $i=1,\ldots, n$, form a
sequence of independent Bernoulli variables. Since
$$\sum_{i=1}^n  \left( p \ln\frac{n}{i}\right)^\ell \sim p^{\ell } n \int_0^1 (\ln 1/t)^{\ell} \d t = p^{\ell } n \ell !,$$
we infer that for $a>0$,
$$p=p_n\sim \frac{a}{n^{1/\ell}}$$
should be the regime in which the largest subpopulations have size
$\ell+1$, and more precisely, then
$$N^{n,p_n}(\ell) \xrightarrow[]{(d)} \poi(\ell!a^{\ell }).$$

These informal calculations are of course far from a rigorous proof;
nonetheless we shall show in the forthcoming Theorem \ref{thm:finitesubp}
that the above weak convergence actually holds. For this, we shall first
develop a general strategy to analyze the asymptotic behavior of the number
of large subpopulations in the next section.

\section{The number of subpopulations exceeding a given size}
\label{sec:strategy}
We will establish limit laws for the number of subpopulations of size
$>x\in\N$ at time $\tau_n$,
\begin{equation}
\label{eq:N-on-E_n-2}
N^{n,p}(x)=\sum_{i=1}^{\infty}\1_{\left\{Y^{(p)}_i(\tau_n)>x\right\}},
\end{equation}
when $p=p_n\rightarrow 0$. The threshold $x$ may be fixed or
$x=x_n\rightarrow\infty$ (depending on the choice of $(p_n)_n$). Roughly
speaking, the key step consists in approximating the population system
$\bY^{(p)}$ at time $\tau_n$ by systems $\bYbelow^{(p)}$ and
$\bYabove^{(p)}$ in which mutants are born at deterministic times, and such
that $\bY^{(p)}$ is bounded by $\bYbelow^{(p)}$ from below and by
$\bYabove^{(p)}$ from above. The advantage of working with deterministic
times comes from the fact that we can in this way decouple populations
sizes from birth times and then deal with independent Bernoulli
variables. In the new systems, we can compute moments of the quantities
corresponding to $N^{n,p}(x)$ more easily. Provided we construct
$\bYbelow^{(p)}$ and $\bYabove^{(p)}$ such that the expected numbers of
subpopulations of a size larger than $x$ match asymptotically in the two
systems, this will ultimately lead to limit statements for $N^{n,p}(x)$ in
the original system.

In fact, it will be sufficient to control $\bY^{(p)}$ on a set of large
probability.  In this regard, recall that convergence in distribution of a
sequence of random variables $(U_n:n\in\N)$ to some random variable $U$ is
equivalent to convergence in distribution of $U_n\1_{E_n}$ to $U$, provided
$(E_n:n\in\N)$ is a sequence of events with $\P(E_n)\rightarrow
1$. Concerning limit results for $N^{n,p}$, we can therefore restrict
ourselves to certain ``good'' events $E_n$ of large probability, which we
specify next.

\begin{lemma}
\label{lem:goodevent}
Assume that $p_n\rightarrow 0$, and let $(k_n)_n\subset \N$ be a sequence
of integers such that $k_n\rightarrow\infty$ and $k_np_n\rightarrow 0$ as
$n\rightarrow\infty$. There exists an event $E_n$ (depending on $n$, $p_n$
and $k_n$) of probability $\P(E_n)\rightarrow 1$, on which the following
holds true:
\begin{enumerate}
\item $\tau_{k_n}=b_{k_n}^{(p)}\leq (\ln k_n)^2$,
\item $(1-k_n^{-1/3})k_n\e^t\,\leq Z(\tau_{k_n}+t)\, \leq
  (1+k_n^{-1/3})k_n\e^t\quad\textup{for all }t\geq 0,$
\item
$
(1-3k_n^{-1/3})k_n\e^t\leq T^{(p_n)}(b_{k_n}^{(p)}+t)\leq
(1+3k_n^{-1/3})k_n\e^t\quad\textup{for all }t\geq 0.
$
\end{enumerate}
\end{lemma}
We call $E_n$ the {\it good event}, and for establishing our limit results
for $N^{n,p_n}(x_n)$, we will work on $E_n$. The proof of the lemma uses
standard estimates on Yule processes and is given in the appendix.

Note that the event $\{\tau_{k_n}=b_{k_n}^{(p)}\}$ under $(a)$ is precisely
the event that the first $k-1$ born individuals (discounting the ancestor
individual) are all mutants. The choice of $k_n$ will depend on our
applications and will be specified later on. Roughly speaking, we choose
$k_n$ in such a way that with high probability, the first $2k_n$ born
individuals do not contribute to $N^{n,p_n}(x_n)$.

Define for $i\in\N$
\begin{align*}
\ob(i)&=\ln i-\ln\left(1-3k_n^{-1/3}\right)-\ln k_n,\\ 
\ub(i)&=\ln i-\ln\left(1+3k_n^{-1/3}\right)-\ln k_n,
\end{align*}
so that
\begin{equation} \label{Eq:beta}
(1-3k_n^{-1/3})k_n\exp\left(\ob(i)\right)=(1+3k_n^{-1/3})k_n\exp\left(\ub(i)\right)=i.
\end{equation} 
Lemma~\ref{lem:goodevent} immediately implies the following control over
birth times on the event $E_n$.
\begin{corollary}
\label{cor:goodevent}
On the good event $E_n$, we have the bounds
$$
\tau_{k_n}+\ub(i)\leq \tau_i\leq \tau_{k_n}+\ob(i).
$$
and
$$
b_{k_n}^{(p)}+\ub(i)\leq b_i^{(p)}\leq b_{k_n}^{(p)}+\ob(i).
$$
\end{corollary}

In the next two sections, we show how to upper and lower bound
$N^{n,p_n}(x_n)$ on the good event $E_n$ for $(x_n)_n$ a sequence of
positive integers (possibly constant). For ease of notation, we write $p$
instead of $p_n$, and $x$ instead of $x_n$. We shall always work on the
event $E_n$, with the properties stated in Lemma~\ref{lem:goodevent}.

\subsection{An upper bound}
\label{sec:upperbound}
Recall the notation \eqref{eq:N-on-E_n-2}. 
We  treat the summands $i$ with $i\leq 2k_n-1$ and $i\geq 2k_n$ separately and
first consider the case $i\geq 2k_n$. From Corollary~\ref{cor:goodevent},
we deduce that for $i\geq 2k_n$,
\begin{equation}
\label{eq:upperbound-tau-n}
\tau_n\leq b_{k_n}^{(p)}+\ob(n)\leq b_i^{(p)}+\ob(n) -\ub(i).
\end{equation}
In the sum~\eqref{eq:N-on-E_n-2}, only the terms with $b_i^{(p)}\leq
\tau_n$ do contribute. Thus we can restrict ourselves to summands with
$2k_n\leq i\leq n^\ast$, where
\begin{equation}
\label{eq:nstar}
n^\ast = \max\left\{i\in\N: \ub(i)\leq \ob(n)\right\}.
\end{equation}
Setting $U_i^{n,p}=Y^{(p)}_i(b_i^{(p)}+\ob(n)-\ub(i))$, we obtain
from~\eqref{eq:upperbound-tau-n} and the fact that $Y^{(p)}_i(t)$ is
monotone increasing in $t$ the almost-sure upper bound
$$\sum_{i=2k_n}^{n^\ast}\1_{\left\{Y^{(p)}_i(\tau_n)>x\right\}}\1_{E_n}\leq \sum_{i=2k_n}^{n^\ast}\1_{\left\{U_i^{n,p}>x\right\}}.
$$
Note that the variables $U_i^{n,p}$, $2k_n\leq i\leq n^\ast$, are independent of each
other, and $U_i^{n,p}$ has the law of the population size of a Yule process
with birth rate $p$ at time $\ob(n)-\ub(i)$ when
started from a single individual. 

The summands with $1\leq i\leq 2k_n-1$ we all bound in the same (rough)
way, by disregarding the individual birth times of the corresponding
subpopulations. More specifically, we remark that on $E_n$,
$$
\tau_n\leq \tau_{k_n}+\ob(n)+\leq\ob(n) +\ln^2(k_n).
$$
so that, defining $U_i^{n,p}=Y^{(p)}_i( b_i^{(p)}+\ob(n)+\ln^2(k_n))$ and
using again monotonicity of $Y_i^{(p)}(t)$, we have almost surely
$$
\sum_{i=1}^{2k_n-1}\1_{\left\{Y^{(p)}_i(\tau_n)>x\right\}}\1_{E_n}\leq \sum_{i=1}^{2k_n-1}\1_{\left\{U_i^{n,p}>x\right\}}.
$$
For $1\leq i\leq 2k_n-1$, the $U_i^{n,p}$'s are independent and identically
distributed according to the law of the population size of a Yule process
with birth rate $p$ at time $\ob(n)+\ln^2(k_n)$ when started
from a single individual. Moreover, the families $(U_i^{n,p}:1\leq i\leq
2k_n-1)$ and $(U_i^{n,p}:2k_n\leq i\leq n^\ast)$ are independent of each
other. Letting
$$
\oS^{n,p}(x)=\sum_{i=1}^{n^\ast}\1_{\left\{U_i^{n,p}>x\right\}},$$
we note that our above estimates yield
\begin{equation}
\label{eq:NleqoSas}
N^{n,p}(x)\1_{E_n}\leq \oS^{n,p}(x)\quad\hbox{almost surely.}
\end{equation}
\subsection{A lower bound}
Our lower bound on $N^{n,p}(x)\1_{E_n}$ begins with the trivial estimate
$$
N^{n,p}(x)\geq \sum_{i=2k_n}^{\infty}\1_{\left\{Y^{(p)}_i(\tau_n)>x\right\}},
$$
that is we will consider only summands with an index $i\geq 2k_n$. From
Corollary~\ref{cor:goodevent}, we see that on the good event $E_n$,
\begin{equation}
\label{eq:lowerbound-tau-n}
\tau_n\geq
\tau_{k_n}+\ub(n)=b_{k_n}^{(p)}+\ub(n)\geq
b_i^{(p)} + \ub(n) - \ob(i).
\end{equation}
Now let
$$
n_\ast = \max\left\{i\in\N: \ob(i)\leq \ub(n)\right\}.
$$
Defining $V_i^{n,p}=Y^{(p)}_i(b_i^{(p)}+\ub(n)-\ob(i))$, we arrive 
with~\eqref{eq:lowerbound-tau-n} at the almost-sure lower bound
$$\sum_{i=2k_n}^{\infty}\1_{\left\{Y^{(p)}_i(\tau_n)>x\right\}}\1_{E_n}\geq \sum_{i=2k_n}^{n_\ast}\1_{\left\{V_i^{n,p}>x\right\}}\1_{E_n}.
$$
The variables $V_i^{n,p}$, $2k_n\leq i\leq n_\ast$, are independent, and
$V_i^{n,p}$ is distributed as the population size of a Yule process with
birth rate $p$ at time $\ub(n)-\ob(i)$ when started from a single
individual. With
$$
\uS^{n,p}(x)=\sum_{i=2k_n}^{n_\ast}\1_{\left\{V_i^{n,p}>x\right\}},$$ we
have shown that
\begin{equation}
\label{eq:NgequSas}
N^{n,p}(x)\1_{E_n}\geq \uS^{n,p}(x)\1_{E_n}\quad\hbox{almost surely.}
\end{equation}
Note that the $V_i^{n,p}$'s are not independent of the event $E_n$. By
construction, we have $V_i^{n,p}\leq U_i^{n,p}$ almost surely for $2k_n\leq
i\leq n_\ast$.

\subsection{Bounds on the mean}
The crucial step of our approach is to control the means
$\E\left[\uS^{n,p}(x)\1_{E_n}\right]$ and
$\E\left[\oS^{n,p}(x)\right]$. In this section, we develop bounds valid
for all regimes. The exact asymptotic analysis will then depend on the
choices of $(p_n)_n$ and $(x_n)_n$ and will be postponed to Section~\ref{sec:limitresults}. 

The behavior of both expectations will be dominated by the following integral.
\begin{lemma}
\label{lem:integral-gamma}
Let $p>0$, $x>0$. Then
$$
\int_0^\infty\left(1-\e^{-ps}\right)^x\e^{-s}\dt s=
\frac{\Gamma(1/p+1)\Gamma(x+1)}{\Gamma(1/p +x +1)},
$$
where $\Gamma$ denotes the Gamma-function.
\end{lemma}
\proof
We let first $u=\e^{-s}$ and then $v=u^p$ to obtain
\begin{align*}
\int_0^\infty\left(1-\e^{-ps}\right)^x\e^{-s}\dt s&=\int_0^1(1-u^p)^x\dt u=
\frac{1}{p}\int_0^1(1-v)^xv^{1/p-1}\dt v\\
&=\frac{1}{p}\textup{B}(1/p,x+1)=\frac{\Gamma(1/p+1)\Gamma(x+1)}{\Gamma(1/p +x +1)},
\end{align*}
where the last two identities follow from well-known expressions for the
Beta- and Gamma-functions B and $\Gamma$, respectively.
\QED

From now on, we let for $n,x\in\N$ and $p>0$,
$$
{\rm I}(n,p,x)=n\frac{\Gamma(1/p+1)\Gamma(x+1)}{\Gamma(1/p +x +1)}.
$$
The following asymptotics as $n\rightarrow\infty$ are immediately derived
from Stirling's formula.
\begin{lemma}
\label{lem:Iasymp}
When $p=p_n\rightarrow 0$ and $x=x_n\rightarrow\infty$, then
$$
{\rm I}(n,p,x)\sim\sqrt{2\pi}n\frac{(1/p)^{1/p+1/2}\,x^{x+1/2}}{(1/p+x)^{1/p+x+1/2}},
$$
whereas if $p=p_n\rightarrow 0$ and $x\in\N$ is fixed, then
$$
{\rm I}(n,p,x)\sim x!np^x.
$$
\end{lemma}
\subsubsection{An upper bound on the mean}
In order to evaluate $\E\left[\oS^{n,p}(x)\right]$, recall the definition
of $U_i^{n,p}$ for $1\leq i\leq n^\ast$.  When evaluated at time $t\geq 0$,
a Yule process with birth rate $p$, started from a single
individual, follows the geometric law with success
probability $\exp(-pt)$. For $1\leq i\leq 2k_n-1$, this implies
$$
\P\left(U_i^{n,p}>x\right)=\left(1-\exp\left(-p(\ob(n)+\ln^2(k_n))\right)\right)^x.
$$
Letting $\Sigma_1^{n,p}(x)=\E[\sum_{i=1}^{2k_n-1}\1_{\{U_i^{n,p}>x\}}]$, we
thus have
\begin{equation}
\label{eq:upperboundmean-1}
\Sigma_1^{n,p}(x)=2k_n\left(1-\exp\left(-p(\ob(n)+\ln^2(k_n))\right)\right)^x.
\end{equation}
In our applications, we will choose $k_n$ such that
$\Sigma_1^{n,p}(x)\rightarrow 0$ as $\rightarrow\infty$, that is, the first
$2k_n-1$ summands are negligible. For the summands with $i\geq 2k_n$, we
get
$$\E\left[\sum_{i=2k_n}^{n^\ast}\1_{\left\{U_i^{n,p}>x\right\}}\right]=\sum_{i=2k_n}^{n^\ast}\left(1-\exp\left(-p(\ob(n)-\ub(i))\right)\right)^x,
$$
and it remains to analyze the sum on the right. We put
$f(n,t)=(1+3k_n^{-1/3})k_n\e^t$ and obtain
\begin{align}
\label{eq:upperboundmean-2}
\sum_{i=2k_n}^{n^\ast}\left(1-\exp\left(-p(\ob(n)-\ub(i))\right)\right)^x&=\sum_{i=2k_n}^{n^\ast}\left(1-\exp\left(-p(\ob(n)-\ub(i))\right)\right)^x\int_{\ub(i-1)}^{\ub(i)}f(n,t)\dt
t\nonumber\\
&\leq
\int_{\ub(2k_n-1)}^{\ub(n^\ast)}\left(1-\exp\left(-p(\ob(n)-t)\right)\right)^xf(n,t)\dt t\nonumber\\
&\leq
\int_0^{\ob(n)}\left(1-\exp\left(-p(\ob(n)-t)\right)\right)^xf(n,t)\dt
t.
\end{align}
In the first step, we used the fact that
$\int_{\ub(i-1)}^{\ub(i)}f(n,t)\dt t=1$ by
\eqref{Eq:beta}, and for the last equality that $\ub(n^\ast)\leq
\ob(n)$ by definition of $n^\ast$, see~\eqref{eq:nstar}. With the change of
variables $s=\ob(n)-t$ and again \eqref{Eq:beta}, the last integral is equal to
$$
\frac{1+3k_n^{-1/3}}{1-3k_n^{-1/3}}\,n\int_0^{\ob(n)}\left(1-\e^{-ps}\right)^x\e^{-s}\dt
  s\leq \left(1+10k_n^{-1/3}\right)n\int_0^{\infty}\left(1-\e^{-ps}\right)^x\e^{-s}\dt
  s,
$$
where the bound on the right holds for $n$ sufficiently large, recalling
that $k_n\rightarrow\infty$. For evaluating the integral on the right side,
we use Lemma~\ref{lem:integral-gamma}. We have shown that for $n$ large enough,
\begin{equation}
\label{eq:N-mean-upperbound}
\boxed{\E\left[N^{n,p}(x)\1_{E_n}\right]\leq \E\left[\oS^{n,p}(x)\right]\leq
  \Sigma_1^{n,p}(x) +
  \left(1+10k_n^{-1/3}\right){\rm I}(n,p,x).}
\end{equation}

\subsubsection{A lower bound on the mean}
We turn to $\E\left[\uS^{n,p}(x)\1_{E_n}\right]$ and write
\begin{equation}
\label{eq:lowerboundmean}
\E\left[\uS^{n,p}(x)\1_{E_n}\right]= \E\left[\uS^{n,p}(x)\right] - \E\left[\uS^{n,p}(x)\1_{(E_n)^c}\right].
\end{equation}
Put $\Sigma_2^{n,p}(x)=\E[\uS^{n,p}(x)\1_{(E_n)^c}].$ In our applications,
$\Sigma_2^{n,p}(x)$ will tend to zero. Indeed, provided $\E[\uS^{n,p}(x)]$
remains uniformly bounded in $n$, the same holds true for its second
moment, see Remark~\ref{rem:strategy} below. We get with Cauchy-Schwarz
$\Sigma_2^{n,p}(x) \leq C(1-\P(E_n))^{1/2} = o(1)$ as $n\rightarrow\infty.$

Concerning the first term on the right side
of~\eqref{eq:lowerboundmean}, we have 
$$\E\left[\uS^{n,p}(x)\right]=\sum_{i=2k_n}^{n_\ast}\left(1-\exp\left(-p(\ub(n)-\ob(i))\right)\right)^x.
$$
Then, with $g(n,t)=(1-3k_n^{-1/3})k_n\e^t$, a similar calculation as
under~\eqref{eq:upperboundmean-2} shows
$$
\sum_{i=2k_n}^{n_\ast}\left(1-\exp\left(-p(\ub(n)-\ob(i))\right)\right)^x
\geq
\int_{\ob(2k_n)}^{\ub(n)}\left(1-\exp\left(-p(\ub(n)-t)\right)\right)^xg(n,t)\dt
t.
$$
Putting $s=\ub(n)-t$, the last integral is bounded from below
by
$$
(1-10k_n^{-1/3})\,n\int_0^{\ub(n)-\ob(2k_n)}\left(1-\e^{-ps}\right)^x\e^{-s}\dt
s.
$$
Now let 
$$
\Sigma_3^{n,p}(x)=\left(1-10k_n^{-1/3}\right)\,n\int_{\ub(n)-\ob(2k_n)}^
\infty\left(1-\e^{-ps}\right)^x\e^{-s}\dt s.
$$
The term $\Sigma_3^{n,p}(x)$ will be negligible in our applications. Using
again Lemma~\ref{lem:integral-gamma}, we have shown that
\begin{equation}
\label{eq:uS-lowerboundmean}
\E\left[\uS^{n,p_n}(x)\right]\geq
\left(1-10k_n^{-1/3}\right){\rm I}(n,p,x)-\Sigma_3^{n,p}(x),
\end{equation}
and therefore
\begin{equation}
\label{eq:N-mean-lowerbound}
\boxed{\E\left[N^{n,p}(x)\1_{E_n}\right]\geq \E\left[\uS^{n,p}(x)\1_{E_n}\right]\geq
  \left(1-10k_n^{-1/3}\right){\rm I}(n,p,x)
  -\left(\Sigma_2^{n,p}(x)+\Sigma_3^{n,p}(x)\right).}
\end{equation}

For the rest of this text, we refer to the quantities $\Sigma_1^{n,p}(x)$
and $\Sigma_2^{n,p}(x)$, $\Sigma_3^{n,p}(x)$ as the error terms of
$\E[\oS^{n,p}(x)]$ and $\E[\uS^{n,p}(x)\1_{E_n}]$, respectively, whereas
the term ${\rm I}(n,p,x)$ involving the Gamma-function (which is the same
for both expectations) is referred to as the main term.
\subsection{General strategy}
Here we outline our general program how to obtain non-degenerate (Poisson)
limit laws for $N^{n,p}(x)$ when $p=p_n\rightarrow 0$ and $x=x_n$ may
depend on $n$ as well. We will control $N^{n,p}(x)$ on the good event $E_n$
in terms of the lower and upper bounds $\uS^{n,p}(x)$ and $\oS^{n,p}(x)$.

Our strategy is based on the following general observation.
\begin{lemma} 
\label{lem:upperconv}
Let $(U_n:n\in\N)$, $(V_n:n\in\N)$ be sequences of $\R$-valued uniformly
  integrable random variables with $U_n\leq V_n$ almost surely for all
  $n\in\N$. Assume $V_n\xrightarrow[]{(d)} V$ for some
  random variable $V$. Then, if 
  $\E[U_n]\rightarrow\E[V]$, we have
  $U_n\xrightarrow[]{(d)} V$ as $n\rightarrow\infty$.
\end{lemma}
\proof Writing $U_n=V_n+(U_n-V_n)$, the claim follows if we show that
$U_n-V_n\rightarrow 0$ in probability. Since $U_n\leq V_n$ almost surely,
we have $\E[|V_n-U_n|]= \E[V_n]-\E[U_n]$. Uniform integrability and the
fact that $V_n\rightarrow V$ in distribution implies $\E[V_n]\rightarrow
\E[V]$, see, e.g.,~\cite[Theorem 3.5]{Bi}. Since $\E[U_n]\rightarrow \E[V]$
by assumption, the proof is complete.  
\QED

We stress that instead of assuming $V_n\xrightarrow[]{(d)} V$ and
$\E[U_n]\rightarrow\E[V]$, one could also {\it assume}
$U_n\xrightarrow[]{(d)} V$ and $\E[V_n]\rightarrow\E[V]$ and then {\it
  deduce} $V_n\xrightarrow[]{(d)} V$ (in fact, for positive integrable
random variables, this can be proved by Fatou's lemma via Skorokhod's
representation theorem, without assuming uniform integrability). However,
for us it is more natural to work under the assumption of
Lemma~\ref{lem:upperconv}.
\begin{remark}
\label{rem:strategy} 
Uniform integrability will not
pose any problem in our setting: If $V_n$ is a sum of indicators of
independent events, then
$$\E\left[V_n^2\right]\leq \E\left[V_n\right] +
\E\left[V_n\right]^2,$$
implying that if $\sup_{n\in\N}\E\left[V_n\right]<\infty$, then
$(V_n:n\in\N)$ has a uniformly bounded second moment as well and is
therefore uniformly integrable.
\end{remark}

We now assume that the sequence $(p_n)_n$ is fixed and satisfies
$p_n\rightarrow 0$. We ask for non-degenerate limits of $N^{n,p_n}(x_n)$.
Recall the definition of $\oS^{n,p_n}(x_n)$. 
Our strategy consists of the following steps. We use the notation from the
preceding sections.  
\newcounter{qcounter}
\begin{list}{{\bf Step \arabic{qcounter}.~}}{\usecounter{qcounter}}
\item Find a sequence $(x_n)_n\subset\N$ such that, for a suitable
  choice of $(k_n)_n$ with $p_nk_n\rightarrow 0$, 
$$
\liminf_{n\rightarrow\infty}\E\left[\uS^{n,p_n}(x_n)\1_{E_n}\right]\geq
\lambda\quad\hbox{and}\quad\limsup_{n\rightarrow\infty}\E\left[\oS^{n,p_n}(x_n)\right]\leq
\lambda$$ for some strictly positive constant $\lambda>0$.  By
construction, this implies
$$
\E\left[\uS^{n,p_n}(x_n)\1_{E_n}\right]\sim\E\left[N^{n,p_n}(x_n)\1_{E_n}\right]\sim\E\left[\oS^{n,p_n}(x_n)\right]\sim
\lambda.
$$
\item Show that for $(k_n)_n$ as under Step $1$, with
$a_n=n^\ast$ and $q_{n,i}=\P\left(U_i^{n,p_n}>x_n\right)$,
  $$\lim_{n\rightarrow\infty}\sum_{i=1}^{a_n}q_{n,i}^2=0.$$
\item Le Cam's inequality applied to the sum $\oS^{n,p_n}(x_n)$ gives
  $\oS^{n,p_n}(x_n)\xlongrightarrow[]{(d)}\poi(\lambda)$. Apply
  Lemma~\ref{lem:upperconv} to deduce that
  $N^{n,p_n}(x_n)\1_{E_n}\xlongrightarrow[]{(d)}\poi(\lambda)$ and hence
  $N^{n,p_n}(x_n)\xlongrightarrow[]{(d)}\poi(\lambda)$ as
  $n\rightarrow\infty$.
\end{list}
Note that the constant $\lambda$ under Step $1$ does not depend on the
sequence $(k_n)_n$. Of course, depending on the point of view, one can also
{\it fix} a sequence of thresholds $(x_n)_n$ and then {\it ask} for a
choice of $(p_n)_n$ such that non-degenerate limits of $N^{n,p_n}(x_n)$
appear.

It is also of interest to understand when $N^{n,p_n}(x_n)\rightarrow 0$ or
$N^{n,p_n}(x_n)\rightarrow\infty$ in probability. In order to prove such
behaviors, we will in the first case apply Markov's inequality and then
show that the expectation converges to zero, while for the second case, we
will make use of the following lemma.
\begin{lemma}
\label{lem:convtoinfinity}
In the setting of Le Cam's inequality, cf. Section~\ref{sec:poisson-heu}, assume that
$\lambda_n=\E[S_n]\rightarrow \infty$ as $n\rightarrow\infty$. Then
$S_n\xrightarrow[]{(p)}\infty$, i.e., for each
$K>0$, $\P(S_n\geq K)\rightarrow 1$ as $n\rightarrow\infty$.
\end{lemma}
The proof follows from an application of the Bienaym\'e-Chebycheff inequality.

We remark that Le Cam's inequality can be used to obtain quantitative
bounds on the rate of convergence. In our case, since we regard
$N^{n,p_n}(x_n)$ only on the good event $E_n$, one would have to optimize the
choice of this event for establishing good bounds. This will not be our concern
here.

\section{Limit results for subpopulation sizes}
\label{sec:limitresults}
We will now work out our strategy explained in the last section. Even though
it can be applied to any choice of $p_n$ such that $p_n\rightarrow
0$, we will restrict ourselves to discussing three regimes of
particular interest. Each of them corresponds to a different limiting behavior of
the ancestral subpopulation, as discussed in Proposition \ref{prop:ancestralsubp}.

\subsection{The regime of bounded subpopulations}
\label{sec:finitesubp}
In Section~\ref{sec:poisson-heu}, we argued heuristically that in the regime
$$
p_n\sim\frac{a}{n^{1/\ell}},\quad a>0\hbox{ and } \ell\in\N \hbox{ fixed,} 
$$
there should be a Poissonian number of subpopulations of size $>\ell$ when
$n\rightarrow\infty$.  We shall now prove this rigorously, together with
the fact that there are no subpopulations of a size strictly larger than
$\ell+1$, and an unbounded number of subpopulations of size $\ell$ (or,
more generally, of size $j$ for each $1\leq j\leq \ell$, see
Corollary~\ref{cor:1}). We point to Theorem 5.4 of~\cite{Bo} for similar
results for percolation on the complete graph, that is, the
Erd\H{o}s-R\'enyi random graph model.

It will be convenient to use the notation
$\Delta^{n,p_n}(\ell)=N^{n,p_n}(\ell-1)-N^{n,p_n}(\ell)$ for the number of
subpopulations of size equal to $\ell$. Recall that $c$ and $C$ denote generic
constants whose values may change from line to line.

\begin{theorem}
\label{thm:finitesubp}
Fix $\ell\in\N$ and $a>0$, and assume $p_n\sim an^{-1/\ell}$. Then,
as $n\rightarrow\infty$,
$$
\Delta^{n,p_n}(\ell+1)\xrightarrow[]{(d)}\poi(\ell!a^\ell).
$$
Moreover, $N^{n,p_n}(\ell+1)\xrightarrow[]{(p)} 0$,
and $\Delta^{n,p_n}(\ell)\xrightarrow[]{(p)} \infty$. 
\end{theorem}

\proof Our first two claims follow if we show that
$N^{n,p_n}(\ell)\xrightarrow[]{(d)}\poi(\ell!a^\ell)$ and
$N^{n,p_n}(\ell+1)\xrightarrow[]{(p)} 0$ as $n$ tends to infinity. We
follow the strategy outlined in Section~\ref{sec:strategy}. We choose
$k_n=\lfloor \ln n\rfloor$ and first bound
$\E\left[\oS^{n,p_n}(\ell)\right]$ from above,
via~\eqref{eq:N-mean-upperbound}. The error term $\Sigma_1^{n,p_n}(\ell)$
is estimated by
$$
\Sigma_1^{n,p_n}(\ell)\leq \frac{3a^\ell\ln^{\ell+1}n}{n}=o(1).
$$
For the main term in~\eqref{eq:N-mean-upperbound}, we obtain with
Lemma~\ref{lem:Iasymp} as $n\rightarrow\infty$,
$$
{\rm I}(n,p_n,\ell)=n\frac{\Gamma(1/p_n+1)\Gamma(\ell+1)}{\Gamma(1/p_n
  +\ell +1)}\sim \ell!n p_n^\ell=\ell!a^\ell+o(1).
$$
The last two display imply
$$\E\left[\oS^{n,p_n}(\ell)\right]\leq \ell!a^\ell
+o(1).$$ We turn to a lower bound on $\E\left[\uS^{n,p_n}(\ell)\1_{E_n}\right]$,
which we will establish via 
expression~\eqref{eq:N-mean-lowerbound}. We first show that the error terms
$\Sigma_2^{n,p_n}(\ell)$ and $\Sigma_3^{n,p_n}(\ell)$ tend to zero.  For
that purpose, recall that by construction,
$$
\E\left[\uS^{n,p_n}(\ell)\right]\leq
\E\left[\oS^{n,p_n}(\ell)\right]\leq C.
$$
Hence the second moment of
$\uS^{n,p_n}(\ell)$ is uniformly bounded as
well, see Remark~\ref{rem:strategy}, so that with Cauchy-Schwarz,
$$
\Sigma_2^{n,p_n}(\ell)\leq C(1-\P(E_n))^{1/2}=o(1).
$$
For the error term $\Sigma_3^{n,p_n}(\ell)$
in~\eqref{eq:N-mean-lowerbound}, we note that
$\ub(n)-\ob(2k_n)\geq \ln n-3\ln\ln n$, whence
\begin{align}
\label{eq:sigma3-error}
\Sigma_3^{n,p_n}(\ell)&\leq n\int_{\ub(n)-\ob(2k_n)}^
\infty\left(1-\e^{-p_ns}\right)^\ell\e^{-s}\dt s\nonumber\\
&\leq \frac{1}{n}+n\int_{\ln n-3\ln\ln n}^
{2\ln n}\left(1-\e^{-p_ns}\right)^\ell\e^{-s}\dt s +
\leq C\frac{\ln^{4\ell}n}{n} = o(1).
\end{align}
Since the main term of~\eqref{eq:N-mean-lowerbound} agrees with that
of~\eqref{eq:N-mean-upperbound}, we get
$
\E\left[\uS^{n,p_n}(\ell)\1_{E_n}\right]\geq \ell!a^\ell - o(1), 
$
and consequently
$$
\E\left[\uS^{n,p_n}(\ell)\1_{E_n}\right]\sim \E[N^{n,p_n}(\ell)\1_{E_n}] \sim
\E\left[\oS^{n,p_n}(\ell)\right] \sim \ell!a^\ell\quad\hbox{ for }n\rightarrow\infty.
$$
Step $1$ of the strategy is therefore established
(with the constant $\lambda$ there given by $\ell!a^\ell$),
and so is Step $2$, since
\begin{align*}
\sum_{i=1}^{n^\ast}\P\left(U_i^{n,p_n}>\ell\right)^2&\leq
\max_{i}\P\left(U_i^{n,p_n}>\ell\right)\E\left[\oS^{n,p_n}(\ell)\right]\\
&\leq C\left(1-\exp\left(-p_n(\ob(n)+\ln^2(k_n))\right)\right)^\ell=o(1).
\end{align*}
We follow Step 3 and obtain
$N^{n,p_n}(\ell)\xrightarrow[]{(d)}\poi(\ell!a^\ell)$ as $n$ tends to
infinity.

For the second part of the theorem, we use Markov's
inequality to obtain
$$\P\left(N^{n,p_n}(\ell+1)\geq 1\right)\leq 
\E\left[N^{n,p_n}(\ell+1)\1_{E_n}\right]+\P(E_n^c)\leq
\E\left[\oS^{n,p_n}(\ell+1)\right] +o(1).
$$
Using again the bound~\eqref{eq:N-mean-upperbound} for
$\E\left[\oS^{n,p_n}(\ell+1)\right]$, we first estimate the error term
$\Sigma_1^{n,p_n}(\ell+1)$ and then the main term similarly to above and
obtain $\E\left[\oS^{n,p_n}(\ell+1)\right]\leq Cn^{-1/\ell}$. This proves
$N^{n,p_n}(\ell+1)\xrightarrow[]{(p)} 0$ as $n\rightarrow\infty$.

It remains to show that $\Delta^{n,p_n}(\ell)\xrightarrow[]{(p)}\infty$,
that is for each $K\in\N$, $\P(\Delta^{n,p_n}(\ell)>K)\rightarrow 1$ as $n$
tends to infinity. Writing
$\Delta^{n,p_n}(\ell)=N^{n,p_n}(\ell-1)-N^{n,p_n}(\ell)$, we have seen that
$N^{n,p_n}(\ell)$ converges in distribution to a Poisson random variable,
so we may prove the claim for $N^{n,p_n}(\ell-1)$ instead of
$\Delta^{n,p_n}(\ell)$. Since $N^{n,p_n}(\ell-1)\1_{E_n}\geq
\underline{S}^{n,p_n}(\ell-1)\1_{E_n}$ almost surely,
cf.~\eqref{eq:NgequSas}, we can estimate
\begin{align*}
\P\left(N^{n,p_n}(\ell-1)\geq K\right)&\geq
\P\left(N^{n,p_n}(\ell-1)\geq K;\,E_n\right)-o(1)\geq
\P\left(\underline{S}^{n,p_n}(\ell-1)\geq K;\,E_n\right)-o(1)\\
&\geq \P\left(\underline{S}^{n,p_n}(\ell-1)\geq K\right)-o(1).
\end{align*}
We analyze the mean of $\underline{S}^{n,p_n}(\ell-1)$
via~\eqref{eq:uS-lowerboundmean} and obtain with Stirling's formula
$$
\E\left[\underline{S}^{n,p_n}(\ell-1)\right] \geq cn^{1/\ell}.
$$
An application of Lemma~\ref{lem:convtoinfinity} shows
$\underline{S}^{n,p_n}(\ell-1)\xrightarrow[]{(p)} \infty$
and thus finishes the proof.\QED

Theorem~\ref{thm:finitesubp} does not tell us how the number of
subpopulations equal to $j$ behave when $j$ is strictly less than
$\ell$. This can however be deduced from the cases $j=\ell$ and $j=\ell+1$,
as we show next.
\begin{corollary}
\label{cor:1}
In the setting of Theorem~\ref{thm:finitesubp}, we have
$\Delta^{n,p_n}(j)\xrightarrow[]{(p)} \infty$ as $n\rightarrow\infty$ for each $1\leq j \leq \ell$.
\end{corollary}
\proof For $j=\ell$, the statement forms already part of the theorem, so we
fix $j\in\N$ with $1\leq j<\ell$. Put $\alpha_n=\lfloor
n^{1-j/\ell}\rfloor$. We will show that the number of
subpopulations of size $j$ that stem from the $\alpha_n$th individual is already 
unbounded as $n\rightarrow\infty$, see Figure 1.

In this regard, we let $Z_{\alpha_n}=(Z_{\alpha_n}(t+\tau_{\alpha_n}):
t\geq 0)$ denote the process that counts the individuals which are
descendants of the $\alpha_n$th individual of $Z$,
i.e. $Z_{\alpha_n}(t+\tau_{\alpha_n})$ is the number of individuals in the
original system at time $t+\tau_{\alpha_n}$ which have the $\alpha_n$th
individual as their common ancestor, no matter whether there are clones or
mutants. It should be clear that $Z_{\alpha_n}$ evolves as a standard Yule
process started from one individual. Next, note that by construction,
$$
\tau_n-\tau_{\alpha_n}\,\,\eqd\,\,\sum_{j=\alpha_n}^{n-1}\frac{1}{j}\mathcal{E}_j
$$
for $(\mathcal{E}_j:j\in\N)$ a sequence of independent standard exponentials.
In particular, 
$$\E\left[\tau_n-\tau_{\alpha_n}\right]= (j/\ell)\ln n +o(1),\quad\textup{Var}\left(\tau_n-\tau_{\alpha_n}\right)=o(1),$$
hence $(\tau_n-\tau_{\alpha_n}) -(j/\ell)\ln n\rightarrow 0$ in probability
as $n\rightarrow\infty$. Since $\e^{-t}Z_{\alpha_n}(t)$ converges almost
surely to a standard exponential variable $\mathcal{E}$ as
$t\rightarrow\infty$, this implies
\begin{equation}
\label{eq:cor1-eq1}
\lim_{n\rightarrow\infty}n^{-j/\ell}Z_{\alpha_n}\left(\tau_n-\tau_{\alpha_n}\right)=
\mathcal{E}\quad\textup{in probability.}
\end{equation}
\begin{figure}[ht]
\centering\parbox{8cm}{\includegraphics[width=7cm]{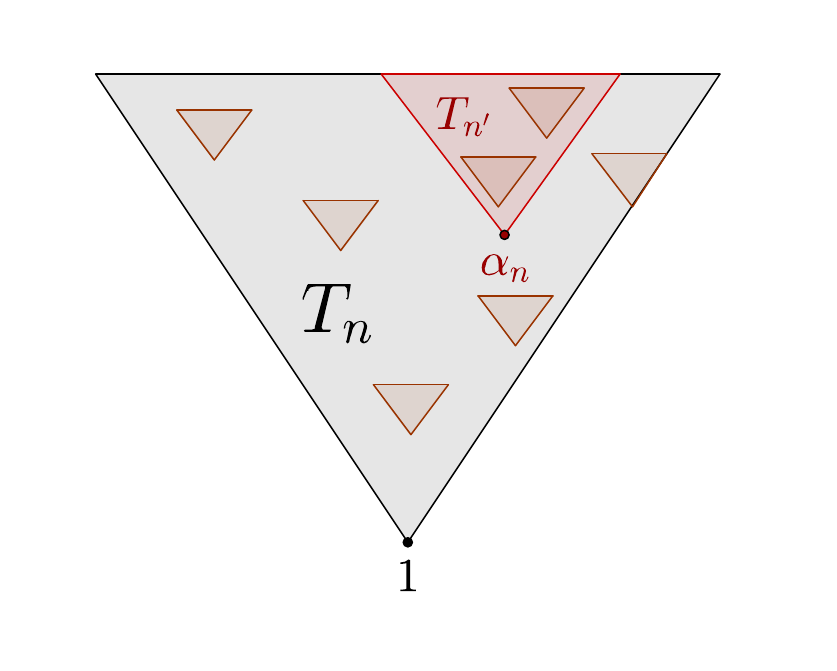}}
\parbox{8cm}{\caption{Schematic of the proof of Corollary~\ref{cor:1}. The
    outer triangle represents the genealogical tree $T_n$ of $Z(\tau_n)$,
    and the small triangles are the subpopulations of size $j$ in
    $\bY^{(p_n)}(\tau_n)$. When $\alpha_n\sim n^{1-j/\ell}$, the subtree
    $T_{n'}$ of $T_n$ rooted at the $\alpha_n$th individual has size
    $n'\approx n^{j/\ell}$. Since $p_n \approx {(n')}^{-1/j}$,
    Theorem~\ref{thm:finitesubp} shows that within $T_{n'}$, the number of
    subpopulations of size $j$ tends to infinity when
    $n\rightarrow\infty$. {\it A fortiori}, the same must hold for
    $\Delta^{n,p_n}(j)$.}}
\end{figure}
Write $N^{n,p_n}_{|\alpha_n}(k)$ for the number of subpopulations of size
$>k\in\N$ in the system $\bY^{(p)}$ stopped at time $\tau_n$, whose common
ancestor is given by the $\alpha_n$th individual, and similarly, define
$\Delta^{n,p_n}_{|\alpha_n}(k)$ by counting only the subpopulations of size
equals $k$ that stem from individual $\alpha_n$. Obviously,
$N^{n,p_n}_{|\alpha_n}(k)\leq N^{n,p_n}(k)$ and
$\Delta^{n,p_n}_{|\alpha_n}(k)\leq \Delta^{n,p_n}(k)$. 
For $K>0$, we estimate
\begin{align}
\label{eq:cor1-eq2}
  \P\left(\Delta^{n,p_n}(j)\geq
    K\right)&\geq\P\left(\Delta^{n,p_n}_{\alpha_n}(j)\geq K\right)
  =\P\left(N^{n,p_n}_{|\alpha_n}(j-1)\geq K +
    N^{n,p_n}_{|\alpha_n}(j)\right)\nonumber\\
  &\geq \P\left(N^{n,p_n}_{|\alpha_n}(j-1)\geq 2 K\right) - \P\left(
    N^{n,p_n}_{|\alpha_n}(j)> K \right).
\end{align}
We will show that the first probability tends to $1$ for each choice of
$K$, while the second can be made as small as we wish provided $K$ is
sufficiently large. 

For that purpose, we remark that conditionally on
$Z_{\alpha_n}\left(\tau_n-\tau_{\alpha_n}\right)=m$, as a consequence of
the dynamics, $N^{n,p_n}_{|\alpha_n}(k)$ for $k\in\N$ has same law as
$N^{m,p_n}(k)$, the number of subpopulations exceeding $k$ in
$\bY^{(p_n)}(\tau_m)$. Moreover, if $m,m'\in \N$ with $m\leq m'$, the
variable $N^{m,p_n}(k)$ is stochastically dominated by $N^{m',p_n}(k)$
(adding individuals can only increase subpopulations). Now fix
$\varepsilon>0$. By~\eqref{eq:cor1-eq1}, we find $n_0\in\N$ and $c_1,
C_1>0$ such that for all $n\geq n_0$, with $m_n=\lfloor
c_1n^{j/\ell}\rfloor$ and $M_n=\lceil C_1n^{j/\ell}\rceil$, the event
$$
A_n=\left\{Z_{\alpha_n}\left(\tau_n-\tau_{\alpha_n}\right)\in
  [m_n,M_n]\right\}$$ has probability at least
$1-\varepsilon$. We first look at the second probability
in~\eqref{eq:cor1-eq2}. By our observations from above, we have for $n\geq n_0$,
$$
\P\left(
    N^{n,p_n}_{|\alpha_n}(j)> K \right)\leq  \P\left(
    N^{n,p_n}_{|\alpha_n}(j)> K\,|\,A_n\right) + \varepsilon\leq
\P\left(
    N^{M_n,p_n}(j)> K \right) + \varepsilon.
$$
Recalling that $p_n\sim an^{-1/\ell}\sim (C_1^{1/j}a)M_n^{-1/j}$, we deduce
from Theorem~\ref{thm:finitesubp} that for $K\in\N$,
$$
\P\left(
    N^{M_n,p_n}(j)> K \right)\rightarrow \P\left(
    \textup{Poi}(C_1j!a^j)>K\right)\quad\textup{as }n\rightarrow\infty.
$$
The right side is smaller than $\varepsilon$ provided $K$ is large
enough. For the first probability in~\eqref{eq:cor1-eq2}, we have similarly
$$
\P\left(N^{n,p_n}_{|\alpha_n}(j-1)\geq 2 K\right) \geq
(1-\varepsilon)\P\left(N^{n,p_n}_{|\alpha_n}(j-1)\geq 2
  K\,|\,A_n\right)\geq (1-\varepsilon)\P\left(N^{m_n,p_n}(j-1)\geq 2 K\right).
$$
By Theorem~\ref{thm:finitesubp}, $N^{m_n,p_n}(j-1)\rightarrow\infty$ in
probability as $n$ tends to infinity, hence the probability on the right
tends to $1$ for each choice of $K$. Since $\varepsilon>0$ was arbitrary,
this concludes the proof of the corollary.  \QED

The next corollary of Theorem~\ref{thm:finitesubp} characterizes the
regime where unbounded subpopulations appear in the limit
$n\rightarrow\infty$. For the sake of clarity, we write $\P_{p_n}$ for the
law of the system $\bY^{(p_n)}$.
\begin{corollary}
\label{cor:2}
Let $(p_n)_n\subset[0,1]$. Then 
$$
\left[\lim_{K\rightarrow\infty}\liminf_{n\rightarrow\infty}\P_{p_n}\left(\exists
    i\in\N\textup{ such that }Y_i^{(p_n)}(\tau_n)>K\right)=1\right] \Leftrightarrow
\left[p_n\,n^{1/\ell}\rightarrow\infty\textup{ for each }\ell\in\N\right].
$$
If one of the statements fails, 
$\lim_{K\rightarrow\infty}\liminf_{n\rightarrow\infty}\P_{p_n}\left(\exists
    i\in\N\textup{ such that }Y_i^{(p_n)}(\tau_n)>K\right)=0.$
\end{corollary}
The proof is a direct application of Theorem~\ref{thm:finitesubp} and left
to the reader.

\subsection{The regime $p_n\sim a\ln^{-1} n$}
Parts $(a)$ and $(b)$ of Proposition~\ref{prop:ancestralsubp} identify $$
p_n\sim\frac{a}{\ln n},\quad a>0\hbox{ fixed,}
$$
as the regime in which the ancestral subpopulation becomes non-trivial. Its
size however remains bounded. What are the sizes of the largest
subpopulations that do appear? As a consequence of 
Theorem~\ref{thm:lognsubp2}, we will see that if we shift
the subpopulation sizes by $-(c_1\ln n+c_2\ln\ln n)$ for some explicit
constants $c_1, c_2>0$, then for any $\ell\in\N$ and any $\varepsilon>0$, we
find $C=C(\ell,\varepsilon)$ such that the $\ell$ largest (shifted) sizes are
contained in $[-C,C]$ with probability at least $1-\varepsilon$, provided
$n$ is large enough.

While Theorem~\ref{thm:lognsubp2} holds true whenever $p_n\sim a\ln^{-1}
n$, we will first prove Theorem~\ref{thm:lognsubp1}, which provides a
stronger result valid for the case $p_n= a\ln^{-1}n$. Here we will compute
a correction $c_3$ of order one to the above shift, such that the number of
subpopulations greater than $y_n=c_1\ln n+c_2\ln\ln n +c_3$ converges to a
Poisson limit along all subsequences $(y_{n(m)})_m$ of $(y_n)_n$, whose fractional
part has a limit as $m$ tends to infinity. Theorem~\ref{thm:lognsubp2} then
readily follows from adapting some estimates used in the proof of Theorem~\ref{thm:lognsubp1}.

Before we give the precise formulation of our results, we need some
preparation.  For $a>0$, define
\begin{equation}
\label{def:fa}
f_a(t) =1+t\ln(at) -\frac{1+at}{a}\ln(1+at),\quad t>0.
\end{equation}
On $(0,\infty)$, $f_a$ is a smooth function with $f_a(t)\rightarrow 1$ as
$t\rightarrow 0$ and $f_a(t)\rightarrow-\infty$ for
$t\rightarrow\infty$. Moreover, since on $(0,\infty)$,
$$
f'_a(t)=\ln(at)-\ln(1+at)<0,
$$
the function $f_a$ is strictly decreasing, and there is a unique
$t^\ast=t^\ast(a)\in(0,\infty)$ for which $f_a(t^\ast)=0$. See Figure 2.

\begin{figure}[ht]
  \centering
  \includegraphics[width=10cm]{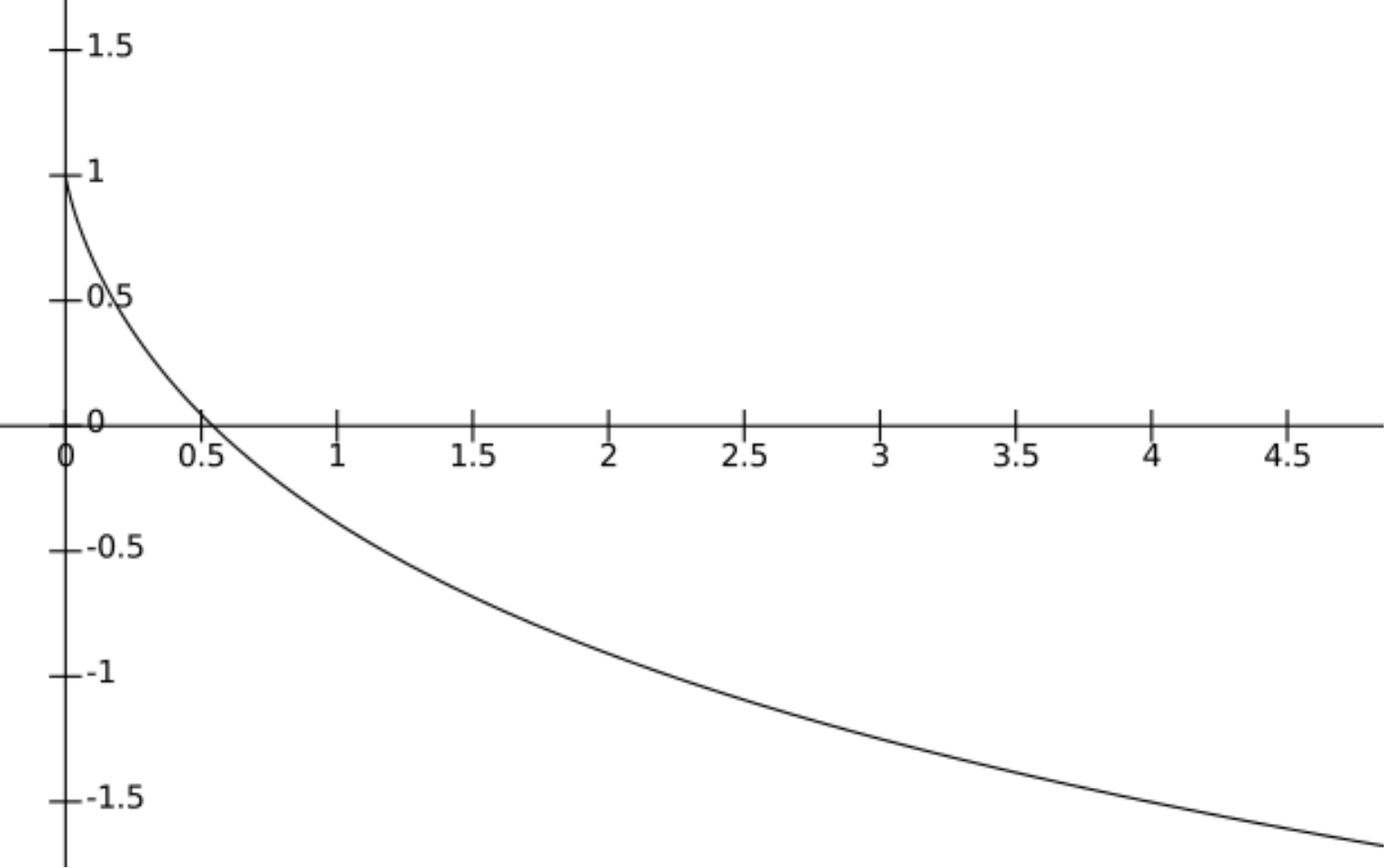}
  \caption{The function $f_a$ for $a=1$.}
 \end{figure}

Now, for $a>0$ and $\lambda>0$ fixed and $t^\ast$ the root of $f_a$, put
\begin{equation}
\label{def:defy_n}
y_n= t^\ast\ln n - \frac{1}{2 f_a'(t^\ast)}\left ( \ln\ln n + \ln\left(\frac{2\pi t^\ast}{\lambda^2(1+at^\ast)}\right)\right) 
\end{equation}

We use the standard notation $\{y\}=y-\lfloor y\rfloor$ to denote the
fractional part of a real $y\in\R$. As we will see in the following
theorem, the barrier $y_n$ defines a precise threshold, in the following
sense: Whenever $(y_{n(m)})_m$ is a subsequence of $(y_{n})_n$ such that
$\lim_{m\rightarrow\infty}\{y_{n(m)}\}=: b\in[0,1)$ exists, then the number
of subpopulations exceeding size $y_{n(m)}$ converges weakly to a
Poisson$(\Lambda(a,b,\lambda))$-distributed random variable with intensity
given by
\begin{equation}
\label{def:Lambda}
\Lambda(a,b,\lambda)=\lambda\left(1+\frac{1}{at^\ast}\right)^b.
\end{equation}
We will see in Corollary~\ref{cor:3} that the restriction to subsequences
with converging fractional part is actually necessary.
\begin{theorem}
\label{thm:lognsubp1}
Assume $p_n= a\ln^{-1} n$ for some $a>0$. Let $\lambda>0$, and define
$y_n=y_n(a,\lambda)$ as in~\eqref{def:defy_n}. Let $(y_{n(m)})_m$ be a
subsequence of $(y_{n})_n$ such that $\{y_{n(m)}\}\rightarrow b$ for some
$b\in[0,1)$. Then, with $\Lambda(a,b,\lambda)$ as above, for
$m\rightarrow\infty$,
$$
N^{n(m),p_{n(m)}}(\lfloor y_{n(m)}\rfloor)\xrightarrow[]{(d)}\poi\left(\Lambda(a,b,\lambda)\right).
$$
\end{theorem}
\proof We again follow the strategy explained in
Section~\ref{sec:strategy}. Recall that we require $p_nk_n\rightarrow 0$,
and we will here choose $k_n=\lfloor \ln^{1/2} n\rfloor$. Putting
$t_n=y_n/\ln n$, the first part of Lemma~\ref{lem:Iasymp} and a small
calculation show that as $n$ tends to infinity,
\begin{equation*} {\rm I}(n,p_n,\lfloor y_n\rfloor)
  \sim \sqrt{\frac{2\pi t_n}{1+at_n}} n\sqrt{\ln n}\frac{(at_n)^{t_n\ln
      n}}{(1+at_n)^{(a^{-1}+t_n)\ln n}}\left(1+\frac{1}{at_n}\right)^{\{t_n\ln n\}}.
\end{equation*}
Taking the logarithm of the right hand side, we arrive at the expression
\begin{equation}
\label{eq:logI}
f_a(t_n)\ln n +\frac{1}{2}\left(\ln\ln n +\ln\left(\frac{2\pi t_n}{1+at_n
      }\right)\right)+\{t_n\ln n\}\ln\left(1+\frac{1}{at_n}\right),
\end{equation}
for $f_a$ as defined under~\eqref{def:fa}. Obviously, $t_n=t^\ast+O(\frac{\ln\ln n}{\ln n})$, and by Taylor's formula,
$$
f_a(t_n)=f_a'(t^\ast)(t_n-t^\ast ) + O\left(\left(\frac{\ln\ln n}{\ln n}\right)^2\right),
$$
so that
$$
f_a(t_n)\ln n = -\frac{1}{2}\left(\ln\ln n +\ln\left(\frac{2\pi t_n}{\lambda^2(1+at_n)}\right)\right)+o(1).
$$
Taking exponentials,~\eqref{eq:logI} and the last display show that for $m$ tending to infinity,
$$
{\rm I}\left(n(m),p_{n(m)},\lfloor y_{n(m)}\rfloor\right) \sim\Lambda(a,b,\lambda).
$$
We next control the error terms $\Sigma_1^{n,p_n}$,
$\Sigma_2^{n,p_n}$ and 
$\Sigma_3^{n,p_n}$. First, recalling that $\ob(n)\leq \ln n +o(1)$,
we have for $n$ sufficiently large,
\begin{align*}
\Sigma_1^{n,p_n}(\lfloor y_n\rfloor)&=2k_n\left(1-\exp\left(-p_n(\ob(n)+\ln^2(k_n))\right)\right)^{\lfloor y_n\rfloor}\nonumber\\
&\leq 2(\ln^{1/2} n) \left(1-\exp(-2a)\right)^{\lfloor y_n\rfloor} = o(1).
\end{align*}
In particular,
$$
\E\left[\oS^{n(m),p_{n(m)}}(\lfloor y_{n(m)}\rfloor)\right]\leq \Lambda(a,b,\lambda) + o(1),
$$
and, with Cauchy-Schwarz and Remark~\ref{rem:strategy},
$\Sigma_2^{n(m),p_{n(m)}}(\lfloor y_{n(m)}\rfloor)=o(1)$. Finally,
\begin{align*}
\Sigma_3^{n,p_n}(\lfloor y_n\rfloor)
&\leq n\int_{\ub(n)-\ob(2k_n)}^
\infty\left(1-\e^{-p_ns}\right)^{\lfloor y_n\rfloor}\e^{-s}\dt s\nonumber\\
&\leq \frac{1}{n}+n\left(1-\e^{-2a}\right)^{\lfloor y_n\rfloor}\int_{\ln n-2\ln\ln n}^
{2\ln n}\e^{-s}\dt s= o(1),
\end{align*}
so that as $m\rightarrow\infty$,
\begin{align*}
\E\left[\uS^{n(m),p_{n(m)}}(\lfloor y_{n(m)}\rfloor)\1_{E_{n(m)}}\right]\sim \E[N^{n(m),p_{n(m)}}(\lfloor y_{n(m)}\rfloor)\1_{E_{n(m)}}] &\sim
\E\left[\oS^{n(m),p_{n(m)}}(\lfloor y_{n(m)}\rfloor)\right]\\ 
&\sim \Lambda(a,b,\lambda).
\end{align*}
The condition under Step $2$ is fulfilled as well, see the estimate for
$\Sigma_1^{n,p_n}$.  Performing Step $3$, it follows that
$N^{n(m),p_{n(m)}}(\lfloor y_{n(m)}\rfloor)$ converges to a
Poisson$(\Lambda(a,b,\lambda))$ random variable. \QED

Theorem~\ref{thm:lognsubp1} implies the following remarkable weak
convergence along subsequences for the size $C_{\ast}^{n,p_n}$ of the
largest subpopulation at time $\tau_n$. For $a>0$, $r\in \R$ and $t^\ast$
the root of $f_a$, put
$$
\lambda_{a,r}=\left(\frac{2\pi t^\ast}{1+at^\ast}\right)^{1/2}\e^{-r/2}.
$$
\begin{corollary}
\label{cor:3}
Let $a>0$, $p_n= a\ln^{-1} n$ and $r\in\R$. Define
$\lambda_{a,r}$ as above, and $y_n=y_n(a,\lambda_{a,r})$ in terms of $\lambda_{a,r}$ as
in~\eqref{def:defy_n}. Let $(y_{n(m)})_m$ be any subsequence of $(y_n)_n$ such that
$\{y_{n(m)}\}\rightarrow b$ for some $b\in[0,1)$. Then, for
$m\rightarrow\infty$, with
$\mu=\ln\left(2\pi/a^{2b}\right)+(1-2b)\ln\left(t^\ast/(1+at^\ast)\right)$,
$$
\P\left(C_{\ast}^{n(m),p_{n(m)}}-t^\ast\ln n(m) + \frac{1}{2f_a'(t^\ast)}\ln\ln
  n(m) \leq
  r\right)\rightarrow \exp\left(-\e^{-(r-\mu)/2}\right).
$$
\end{corollary}
\proof The probability on the left
hand side is equal to
$$
\P\left(C_{\ast}^{n(m),p_{n(m)}}\leq \lfloor y_{n(m)}\rfloor\right)
=\P\left(N^{n(m),p_{n(m)}}(\lfloor
  y_{n(m)}\rfloor)=0\right)\sim\P\left(\poi(\Lambda(a,b,\lambda_{a,r}))=0\right),
$$
with $\Lambda(a,b,\lambda_{a,r})$ defined as in~\eqref{def:Lambda}. The
last asymptotics holds thanks to Theorem~\ref{thm:lognsubp1}. The claim
follows.  \QED

\begin{remark}
  The limit expression on the right hand side in the above corollary is the
  distribution function of the Gumbel law with location parameter $\mu$ and
  scale parameter $2$. In view of our heuristics explained in
  Section~\ref{sec:poisson-heu} and of what is known about the maximum of
  $n$ independent geometrically distributed random variables, it should not
  come as a surprise that a Gumbel distribution appears in the limit for
  the recentered size of the largest subpopulation. Since $C_\ast^{n,p_n}$
  is a discrete random variable and the Gumbel law is a continuous
  distribution, convergence along the full sequence cannot hold.
\end{remark}

It is instructive to compare our previous two results with Theorem 5.10
and Corollary 5.11 of~\cite{Bo} for the Erd\H{o}s-R\'enyi random graph model.

From the proof of Theorem~\ref{thm:lognsubp1}, we easily derive that in the
general case $p_n\sim a\ln^{-1} n$, the sizes of the largest subpopulations
are concentrated around
\begin{equation}
\label{def:defb_n}
b_n= t^\ast\ln n - \frac{1}{2 f_a'(t^\ast)}\ln\ln n.
\end{equation}

\begin{theorem}
\label{thm:lognsubp2}
Assume $p_n\sim a\ln^{-1} n$ for some $a>0$. Define $b_n=b_n(a)$ as
in~\eqref{def:defb_n}, and let $(x_n)_n$ be a sequence of
positive integers. Then the following holds as $n\rightarrow\infty$.
\begin{enumerate}
\item If $(x_n-b_n)\gg 1$, then
  $
N^{n,p_n}(x_n)\rightarrow 0.
$
\item If $(b_n-x_n)\gg 1$, then
$
N^{n,p_n}(x_n)\rightarrow \infty.$
\end{enumerate}
\end{theorem}
\proof We only have to adapt the estimates given in the proof of
Theorem~\ref{thm:lognsubp1}, so we will only sketch the necessary
modifications. We again choose $k_n=\lfloor \ln^{1/2} n\rfloor$.  We 
treat $(a)$ and $(b)$ together, and in this regard, we note that for $(b)$, by
monotonicity of $N^{n,p_n}(x)$ we can assume that $x_n\gg \ln\ln n$.  With $y_n$ everywhere
replaced by $x_n$, we deduce from~\eqref{eq:logI}
that if $(x_n-b_n)\gg 1$, then ${\rm I}(n,p_n,x_n)\rightarrow -\infty$,
whereas if $(b_n-x_n)\gg 1$, then ${\rm I}(n,p_n,x_n)\rightarrow \infty$. The
error term $\Sigma_1^{n,p_n}(x_n)$ is seen to be of order $o(1)$ with
exactly the same argument as in the proof of Theorem~\ref{thm:lognsubp1},
and so is the error term $\Sigma_3^{n,p_n}(x_n)$, using here that $x_n\gg
\ln\ln n$. 

Now if $(x_n-b_n)\gg 1$, we have by Markov's inequality and~\eqref{eq:N-mean-upperbound},
$$
\P\left(N^{n,p_n}(x_n)\geq 1\right)\leq \E\left[\oS^{n,p_n}(x_n)\right] +
o(1) = o(1),
$$
and if $(b_n-x_n)\gg 1$, then, as in the second part of the proof of
Theorem~\ref{thm:finitesubp}, for any $K\in\N$,
$$
\P\left(N^{n,p_n}(x_n)\geq K\right)\geq \P\left(\uS^{n,p_n}(x_n)\geq K\right)-o(1).
$$
The probability on the right tends to $1$ by an appeal to Lemma~\ref{lem:convtoinfinity}.
\QED
\subsection{The regime $\ln^{-1} n\ll p_n\ll 1$.}
In the regime $p_n\gg 1/\ln n$, the ancestral subpopulation grows like
$n^{p_n}$, see part $(c)$ of Proposition~\ref{prop:ancestralsubp}. As we
shall prove in the following theorem, when $1/\ln n\ll p_n\ll 1$, the sizes
of the largest subpopulations are to first order given by
$\e^{-1}p_n^{-1}n^{p_n}$. For the case $\ln\ln n/\ln n\leq p_n\ll 1$ and
$\lambda>0$, we will show that the number of subpopulations greater than
\begin{equation}
\label{def:defx_n}
x_n =x_n(p_n,\lambda)=
\Bigl\lfloor{\left(\lambda^{-1}\sqrt{2\pi p_n^{-1}}\right)}^{p_n}\e^{-1}p^{-1}_nn^{p_n}\Bigr\rfloor
\end{equation} is asymptotically Poisson$(\lambda)$-distributed. Note that
as $n\rightarrow\infty$, $g_\lambda(p_n)=(\lambda^{-1}\sqrt{2\pi
  p_n^{-1}})^{p_n}\sim 1$, cf. Figure 3.

\begin{figure}[ht]
\centering\parbox{8cm}{\includegraphics[width=7.5cm]{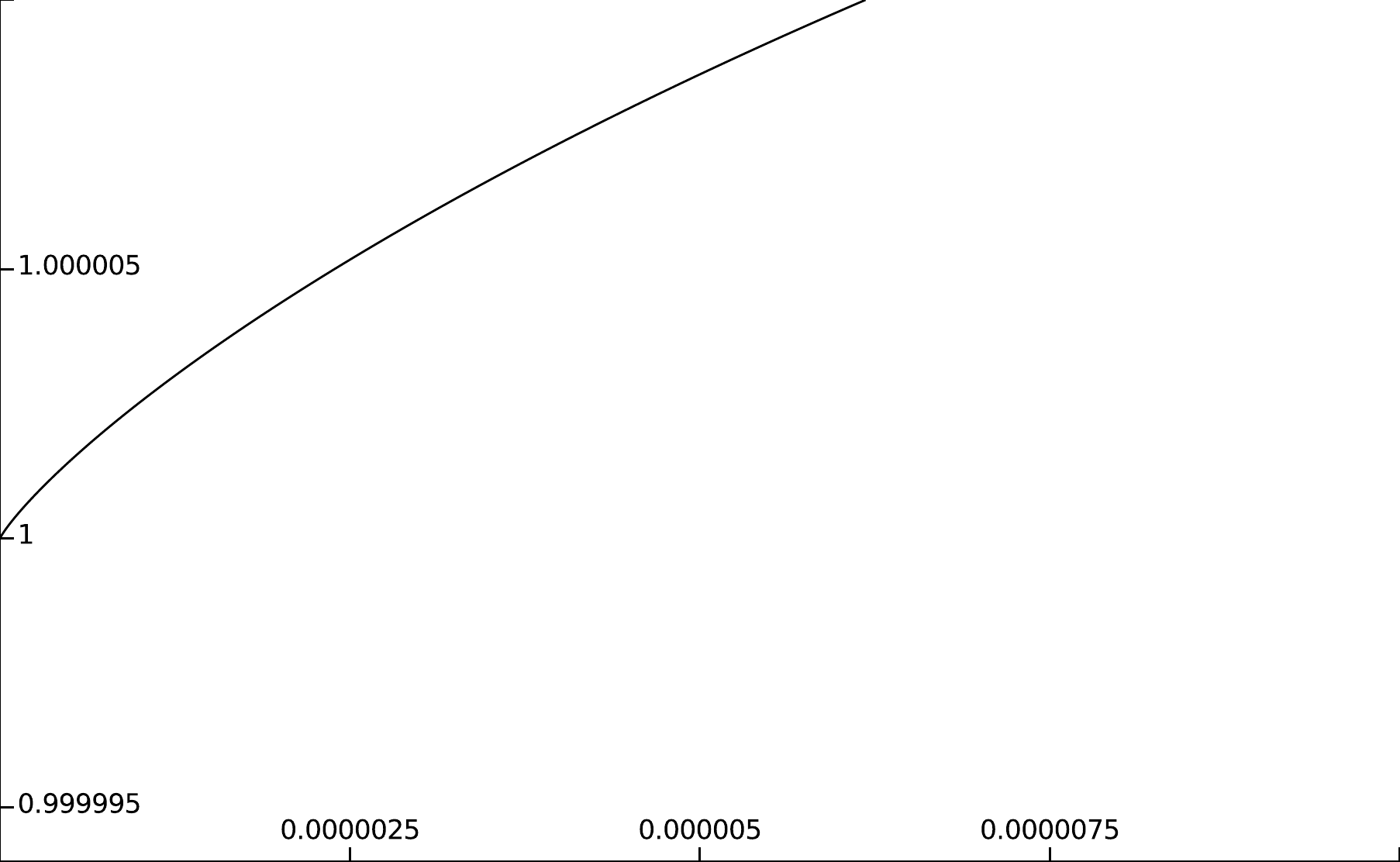}}
\parbox{8cm}{\includegraphics[width=7.5cm]{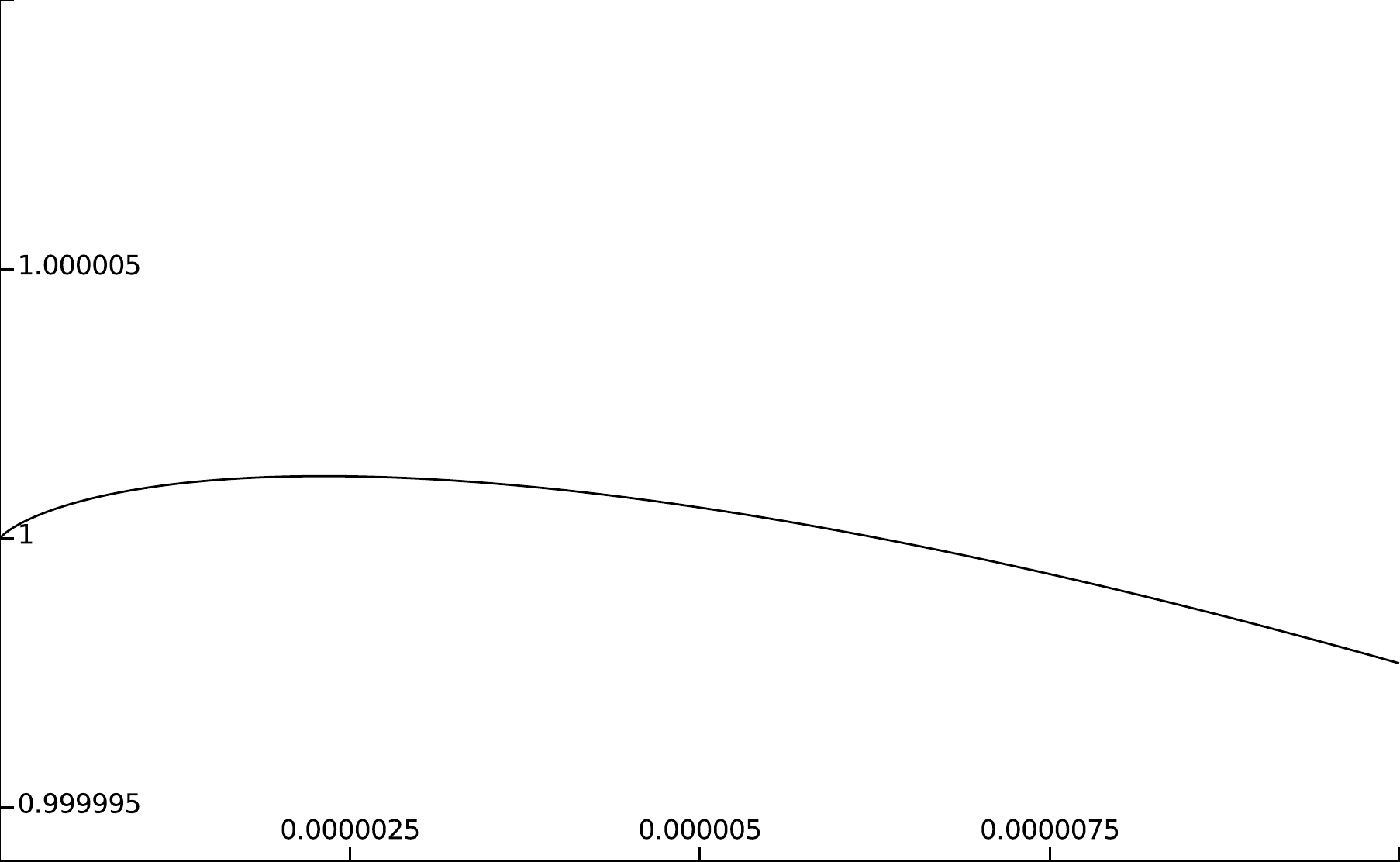}}
\caption{The function $g_\lambda(p)$ on $[0,10^{-5}]$ for $\lambda=200$ (left side)
  and $\lambda=1000$ (right side).}
\end{figure}

As it will become clear from the proof, we require $p_n\geq \ln\ln n/\ln n$
(and not merely $p_n\gg 1/\ln n$) only to simplify the calculation; without
this assumption, the exact behavior of $p_n$ has to be taken into account
more carefully.
\begin{theorem}
\label{thm:greatsubp}
Assume $\ln^{-1}n\ll p_n\ll 1$. Then the following holds as
$n\rightarrow\infty$.
\begin{enumerate}
\item If $u>\e^{-1}$, then $N^{n,p_n}(\lfloor
  up_n^{-1}n^{p_n}\rfloor)\xrightarrow[]{(p)} 0$. If $u<\e^{-1}$, then 
  $N^{n,p_n}(\lfloor up_n^{-1}n^{p_n}\rfloor)\xrightarrow[]{(p)} \infty$.
\item Assume additionally $p_n\geq \ln\ln n/\ln n$ for large $n$. Then,
  with $x_n$ as in~\eqref{def:defx_n},
 $$
N^{n,p_n}(x_n)\xrightarrow[]{(d)}\poi(\lambda).
$$
\end{enumerate}
\end{theorem}
\proof We work with the choice $k_n=\lfloor p_n^{-1/2}\rfloor$ and follow
the steps of the strategy presented in Section~\ref{sec:strategy}.  Since
our proof is similar to the proofs of Theorem~\ref{thm:finitesubp} and 
Theorem~\ref{thm:lognsubp1}, we do not provide every detail.

We fix a real number $u>0$ and let $y_n= \lfloor u p_n^{-1}n^{p_n}\rfloor$,
$n\in\N$.  For the main term ${\rm I}(n,p_n,y_n)$, Lemma~\ref{lem:Iasymp}
shows
\begin{equation}
\label{eq:thm3-mainterm}
{\rm I}(n,p_n,y_n)\sim
\sqrt{2\pi}\,n\frac{(1/p_n)^{1/p_n+1/2}\, y_n^{y_n+1/2}}{(1/p_n+y_n)^{1/p_n+y_n+1/2}}.
\end{equation}
Note that since $p_n\gg 1/\ln n$, we have $p_ny_n\rightarrow\infty$
as $n\rightarrow\infty$. By Taylor's formula,

$$
\ln\left(\frac{1}{p_n}+y_n\right) =
\ln(y_n)+\frac{1}{p_ny_n}+O\left(\frac{1}{p^2_ny^2_n}\right).
$$
From the last display, we deduce after a short calculation that the logarithm
of the right hand side of~\eqref{eq:thm3-mainterm}
behaves asymptotically like
\begin{align}
\label{eq:thm3-logmainterm}
\lefteqn{\frac{1}{2}\ln(2\pi) +\ln n - \left(\frac{1}{p_n}+\frac{1}{2}\right)\ln p_n - \frac{1}{p_n}\ln(y_n) -\frac{1}{p_n} +
O\left(\frac{1}{p^2_ny_n}\right)}\nonumber\\
&\quad\quad\quad\quad= -\frac{1}{p_n}\left(1+\ln(u)+O(n^{-p_n})\right)-\frac{1}{2}\ln p_n+\frac{1}{2}\ln(2\pi). 
\end{align}
In particular, if $u>\e^{-1}$, then the right hand side diverges to
$-\infty$, whereas if $u<\e^{-1}$, it diverges to $+\infty$. Consequently,
${\rm I}(n,p_n,y_n)\rightarrow 0$ in the first and ${\rm
  I}(n,p_n,y_n)\rightarrow \infty$ in the second case. In order to finish
the proof of $(a)$, it remains to convince ourselves that the error terms
$\Sigma_1^{n,p_n}(y_n)$ and $\Sigma_3^{n,p_n}(y_n)$ do not contribute when
$n$ tends to infinity.  For the first one, we have
$$
\Sigma_1^{n,p_n}(y_n)\leq
2p_n^{-1/2}\left(1-\exp\left(-p_n(\ob(n)+\ln^2(\lfloor p_n^{-1/2}\rfloor))\right)\right)^{y_n}.
$$
Since $\ob(n)\leq \ln n + o(1)$, we see from taking logarithms that
$\Sigma_1^{n,p_n}(y_n)=o(1)$, as wanted. The error term
$\Sigma_3^{n,p_n}(y_n)$ is readily seen to be of order $o(1)$ as well, and
the claims under $(a)$ now follow from the same arguments as in the proof
of Theorem~\ref{thm:lognsubp2} (or Theorem~\ref{thm:finitesubp}).

We turn to $(b)$ and fix $\lambda>0$. If $p_n\geq \ln\ln n/\ln n$ for $n$
large, then, with $x_n$ defined as in~\eqref{def:defx_n},
$$
O\left(\frac{1}{p_n n^{p_n}}\right)=o(1).
$$
Performing the above calculation with $y_n$ replaced by $x_n$ and $u$
by 
$$
u_n= {\left(\lambda^{-1}\sqrt{2\pi p_n^{-1}}\right)}^{p_n}\e^{-1},
$$
we deduce from~\eqref{eq:thm3-logmainterm} that the logarithm of right side
in~\eqref{eq:thm3-mainterm} behaves as $\ln \lambda +o(1)$, and thus
$$
{\rm I}(n,p_n,x_n)\sim\lambda.
$$
Part $(b)$ now follows from the same line of reasoning as in the proof of
Theorem~\ref{thm:lognsubp1}, using that all the error terms
$\Sigma_i^{n,p_n}(x_n)$, $i=1,2,3$, are of order $o(1)$.\QED

\section{Subcritical percolation on random recursive trees}
\label{sec:app-percolation}
In this final section, we make the precise link between the population
system ${\bf Y}^{(p)}$ defined in Section~\ref{sec:Yule} and percolation on
random recursive trees.

A recursive tree with vertex set $[n]=\{1,\ldots,n\}$ is a tree rooted at
$1$ with the property that the labels along the unique path from the root
to any other vertex form an increasing sequence. A tree chosen uniformly at
random among all these $(n-1)!$ recursive trees is called random recursive
tree and denoted $T_n$.  The study of Bernoulli bond percolation on large
random recursive trees can be traced back to the analysis of an algorithm
for isolating the root by Meir and Moon \cite{MM}; see also Drmota {\it et
  al.} \cite{DIMR} , Iksanov and M\"ohle \cite{IM} and Kuba and Panholzer
\cite{KP} for more recent developments in this direction. It further plays
a key role in the construction of the Bolthausen-Sznitman coalescent by
Goldschmidt and Martin \cite{GM}.

We choose $p=p_n\in(0,1)$ and write $C_1^{n,p},C_2^{n,p},\ldots$ for the
percolation clusters of a Bernoulli bond percolation on $T_n$ with
parameter $p$.  That is, we erase each edge of $T_n$ with probability $1-p$
and independently of each other, and enumerate the connected components
according to the increasing order of the label of their root vertex (i.e.
their vertex with the smallest label).

We use the convention that $C_i^{n,p}=\emptyset$ if there are less than $i$
connected components after percolation on $T_n$. With our ordering,
$C_1^{n,p}$ always represents the root cluster containing $1$; $C_2^{n,p}$
is the cluster rooted at the smallest vertex which does not belong to the
root cluster $C_1^{n,p}$, and so on.  We write $|C_i^{n,p}|$ for the number
of vertices of $C_i^{n,p}$. Then the following holds.

\begin{lemma}
\label{lem:equivRRTYule}
Let $p\in [0,1]$.
The families $(|C_{1}^{n,p}|,|C_{2}^{n,p}|,\ldots)$ and
$(Y_1^{(p)}(\tau_n),Y_2^{(p)}(\tau_n),\ldots)$ have the same law.
\end{lemma}
\proof We construct percolation on a random recursive tree as a growth
process in continuous time as follows. At time $t=0$, we start from the
singleton $\{1\}$. Given percolation with parameter $p$ on a random
recursive tree on $[k]$, $k\geq 1$, has been constructed, we equip each
vertex $i\in [k]$ with an independent exponential clock $\mathcal{E}_i$ of
parameter $1$. At time $\min_{i=1,\ldots,k}\mathcal{E}_i$, we flip a coin
with heads probability $p$. If head shows up, we attach a vertex labeled
$k+1$ to the vertex with label
arg$\min_{i=1,\ldots,k}\mathcal{E}_i$. Otherwise, we add vertex $k+1$ to
the system without connecting it to any other vertex. It follows from the
construction of random recursive trees and the independence of the coin
flips that we observe at the instant when the $n$th vertex is incorporated
a Bernoulli bond percolation on $T_n$ with parameter $p$. Moreover, if we
keep track of the sizes
$\bC^{(p)}(t)=(|C_1^{(p)}(t)|,|C_2^{(p)}(t)|,\ldots)$ of the growing
percolation clusters, where $|C_i^{(p)}(t)|$ stores the size of the $i$th
cluster at time $t$ and clusters are ordered according to their birth
times, we obtain a Markov chain with initial state
$\bC^{(p)}(0)=(1,0,\ldots)$ and transition rates at time $t\geq 0$ for
$c_1,\ldots,c_k\in\N$ given by
\begin{align*}
  (c_1,\ldots,c_k,0,\ldots) &\mapsto (c_1,\ldots,c_k,1,0,\ldots)\quad
  \textup{at rate }(1-p)(c_1+\ldots+ c_k),\\
  (c_1,\ldots,c_k,0,\ldots)&\mapsto(c_1,\ldots,c_{i-1},c_i+1,c_{i+1},\ldots,c_k,0,\ldots)
  \quad \textup{at rate }pc_i.
\end{align*}
From the very construction of $\bY^{(p)}$, it follows that the processes
$(\bY^{(p)}(t):t\geq 0)$ and $(\bC^{(p)}(t):t\geq 0)$ grow according to the
same dynamics, and the statement follows.
\QED

If we think of the individuals of $\bY^{(p)}$ as being labeled
$1,2,\ldots$, according to the increasing order of their birth times,
there is in fact a one-to-one correspondence
between subpopulations and clusters that respects the genealogy, as
illustrated by Figure 4.

\begin{figure}[ht]
\centering\parbox{8cm}{\includegraphics[width=7cm]{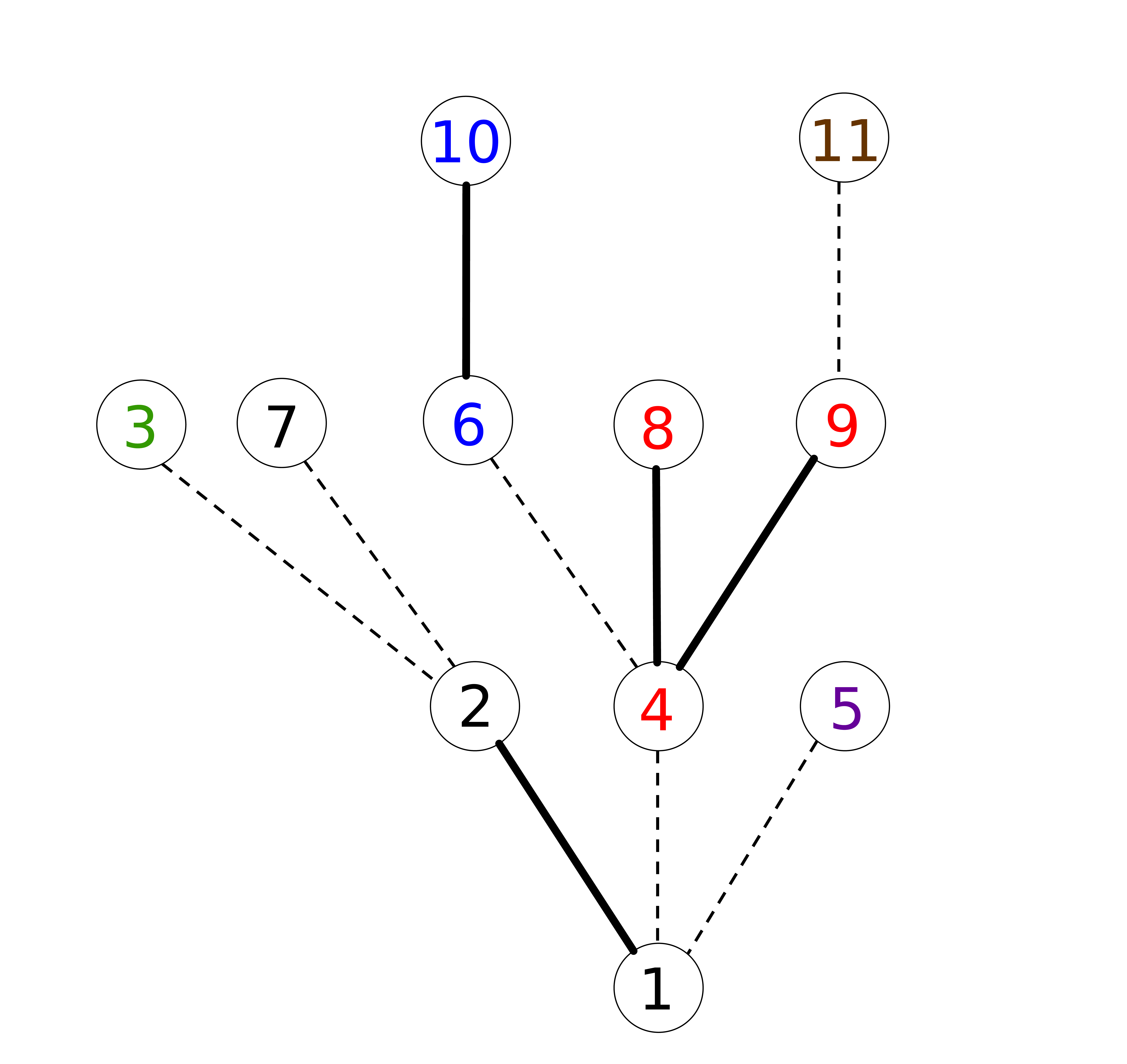}}
\parbox{8cm}
 {\caption{The genealogical tree of a Yule process with mutation, stopped
    after the birth of the $11$th individual. Subpopulations of different
    genetic types have different colors, and the edges connecting mutants
    to their parent are represented by dashed lines.  Alternatively, the
    subpopulations can be viewed as the clusters of a Bernoulli bond
    percolation on a recursive tree on $\{1,\ldots,11\}$.}}
\end{figure}

With the above lemma at hand, the results from
Section~\ref{sec:limitresults} have a direct interpretation in terms of
clusters stemming from a percolation on $T_n$ with parameter
$p_n$. For example, Corollary~\ref{cor:2} gives a necessary and sufficient
condition on the percolation parameter $p_n$ such that in the limit
$n\rightarrow\infty$, clusters of unbounded size appear.

More generally, the strategy developed in Section~\ref{sec:strategy} gives
a concise tool to decide whether for a given percolation sequence
$p_n\rightarrow 0$ and thresholds $x_n$, there are clusters of a size of
order $x_n$, and if so, how many. Indeed, as exemplified in
Section~\ref{sec:limitresults}, basically one only has to check the
asymptotic behavior of the expression ${\rm I}(n,p_n,x_n)$ introduced below
Lemma~\ref{lem:integral-gamma}.

This work completes the study of percolation on random recursive
trees initiated in~\cite{Be2}, where percolation with supercritical
parameter $p_n\sim 1-a/\ln n$, $a>0$ fixed, is studied. In~\cite{Be4},
non-Gaussian fluctuations of the root cluster were proved, and the analysis
of percolation clusters was extended in~\cite{Ba} to all regimes
$p_n\rightarrow 1$. The works~\cite{Be2,Ba} do additionally contain
information on the genealogy of clusters (and not merely on their sizes), a
type of question we did not investigate here.

As a common feature when $p_n\rightarrow 1$, the root cluster containing
$1$ has size $\sim n^{p_n}$, while the sizes of the next largest clusters
are of order $(1-p_n)n^{p_n}$. This motivated a sub-classification of
regions into weakly supercritical [$1/\ln n\ll 1-p_n\ll 1$] when
$|C_{1}^{n,p_n}|\ll n$, supercritical [$1-p_n\sim t/\ln n$, $t>0$ fixed]
when $|C_{1}^{n,p_n}|\sim\e^{-t} n$, and strongly supercritical [$1-p_n\ll
1/\ln n$] when $|C_{1}^{n,p_n}|\sim n$.

The regime of constant $p\in(0,1)$ may best be termed ``critical'', since
in this case, the root cluster loses its dominating role (with respect to
the size). Some results can be found in~\cite{Mo} and~\cite{BaBe}, although
the motivation there is somehow different. For similar reasons, we term the
regime $p=p_n\rightarrow 0$ considered here ``subcritical''. The
classification into different regimes draws on the terminology which is
usually used to describe the Erd\H{o}s-R\'enyi random graph model (see,
e.g.,~\cite{AlSp}).

We stress that there are other natural families of trees which can
be grown according to a probabilistic evolution algorithm, like,
e.g., scale-free random trees or $b$-ary random increasing trees. Indeed,
branching systems with mutations were used in~\cite{BeBr} to study
percolation  on scale-free random trees when $p_n\sim a/\ln n$, and
in a similar fashion by Berzunza in~\cite{Ber} for $b$-ary
random increasing trees. In the case of random recursive trees, the
underlying population system is particularly simple, so we restricted our
discussion of the subcritical regime to these trees, but we certainly
expect similar results to hold true for other classes of increasing tree
families.

\begin{appendix}
\section{Proof of Lemma~\ref{lem:goodevent}}
\label{sec:appendix}
In this appendix, we will  construct an event $E_n$ with the ``good''
properties specified in Lemma~\ref{lem:goodevent}. We first collect some
properties of Yule processes.

\subsection{Some properties of Yule processes}
For $p\in [0,1]$, denote by $Y^{(p)}=(Y^{(p)}(t):t\geq 0)$ a Yule process
with birth rate $p$ per unit population size. Given $y\in\N$, write $\P_y$
for its law under which $\P_y(Y^{(p)}(0)=y)=1$, and similarly $\E_y$ for
its expectation. It is well-known that $W^{(p)}(t)=\e^{-pt}Y^{(p)}(t)$,
$t\geq 0$, is a square-integrable martingale. Under $\P_y$ and for $p>0$,
its terminal value $W^{(p)}(\infty)$ is distributed as a sum of $y$
standard exponentials and follows thus the Gamma$(y,1)$-law. Note that
Doob's inequality applied to the square-integrable martingale
$W^{(p)}(t)-y$ gives

\begin{equation}
\label{eq:Yule-distancetostart}
\E_y\left[\sup_{t\geq 0}\big|W^{(p)}(t)-W^{(p)}(0)\big|^2\right]\leq
4\E_y\left[|W^{(p)}(\infty)-W^{(p)}(0)|^ 2\right]= 4y.
\end{equation}

We now work in the setting of Section~\ref{sec:Yule}. Recall the
definition of the process $T^{(p)}=(T^{(p)}(t):t\geq 0)$ counting the
number of different genetic types in the population system $\bY^{(p)}$,
starting from $T^{(p)}(0)=1$. For $i\in\N$, $\tau_{i} = \inf\left\{t\geq 0:
  Z(t)= i\right\}$ denotes the birth time of the $i$th individual in the
system $\bY^{(p)}$ (counting both clones and mutants), with $Z$ denoting
the underlying standard Yule process. The time $b_i^{(p)}$ denotes the
first time when an individual of genetic type $i$ appears, i.e., the
$b_i^{(p)}$'s are the jump times of $T^{(p)}(t)$. We need the following
control over $T^{(p)}$.
\begin{lemma}
\label{lem:numberoftypes}
Let $p\in(0,1)$. For $k\in\N$, let $\P_{|k}$ be the conditional law given
$\{\tau_k=b_k^{(p)}\}$ (i.e. the first $k-1$ children are all mutants), and let $\E_{|k}$ be its expectation. For every $t\geq 0$, we have
$$
\E_{|k}\left[\sup_{0\leq s\leq
    t}\left|T^{(p)}(\tau_k+s)-\left(k+\int_{\tau_k}^{\tau_k+s}(1-p)Z(r)\dt
      r\right)\right|^2\right]\leq 4(1-p)k(\e^{t}-1),
$$
and
$$
\E_{|k}\left[\sup_{s\geq
    t}\e^{-5s/3}\left|T^{(p)}(\tau_k+s)-\left(k+\int_{\tau_k}^{\tau_k+s}(1-p)Z(r)\dt
      r\right)\right|^2\right]\leq 324(1-p)k\,\e^{-2t/3}.
$$
\end{lemma}

\proof Let $\mathcal{F}_t$ denote the natural filtration generated by the
process $(\bY^{(p)}(t):t\geq 0)$. Both processes $T^{(p)}$ and $Z$ are
$\mathcal{F}_t$-adapted, and $\tau_k$ is a stopping time. We work now under
the probability measure $\P_{|k}$. From the strong
Markov property and the dynamics of $\bY^{(p)}$ described above, we see that
$M^{(p)}(s)=T^{(p)}(\tau_k+s)-(k+\int_{\tau_k}^{\tau_k+s}(1-p)Z(r)\dt r)$
is a martingale with $M^{(p)}(0)=0$. Its quadratic variation is given by
$[M^{(p)}](s)=T^{(p)}(\tau_k+s)-k$, and its second moment takes the form
\begin{align*}
  \E_{|k}\left[|M^{(p)}(s)|^2\right]&=\E_{|k}\left[[M^{(p)}](s)\right]=\E_{|k}\left[T^{(p)}(\tau_k+s)\right]-k\\
  &=\E_{|k}\left[\int_{\tau_k}^{\tau_k+s}(1-p)Z(r)\dt
    r\right]=\int_{0}^{s}(1-p)\E[Z(\tau_k+r)]\dt r= (1-p)k(\e^{t}-1).
\end{align*}
For the last equality, we used the strong Markov property, which entails
that $(Z(\tau_k+r):r\geq 0)$ is independent of $\mathcal{F}_{\tau_k}$ and
distributed as a standard Yule process started from $k$ individuals. An
appeal to Doob's inequality gives
$$
\E_{|k}\left[\sup_{0\leq s\leq t}|M^{(p)}(s)|^2\right] \leq 4 \E_{|k}\left[|M^{(p)}(t)|^2\right]=4(1-p)k(\e^{t}-1).
$$
For the second statement, we bound for every $n\in\N$
$$
\sup_{n\leq s< n+1}\e^{-5s/6}|M^{(p)}(s)|\leq \e^{-5n/6}\sup_{n\leq s< n+1}|M^{(p)}(s)|.
$$
By Doob's inequality, see the first part, the $L^2$-norm of the right side
(with respect to the conditional measure $\P_{|k}$)
is bounded by
$$
\e^{-5n/6}\left\|\sup_{n\leq s< n+1}|M^{(p)}(s)|\right\|_2 \leq
2\e^{-5n/6}\sqrt{(1-p)k\e^{n+1}}\leq 2\e^{1/2}\sqrt{(1-p)k}\,\e^{-n/3}.
$$
Thus, summing over $n$ and applying the triangle inequality,
$$\left\|\sup_{s\geq t}\e^{-5s/6}|M^{(p)}(s)|\right\|_2\leq
2\e^{1/2}\sqrt{(1-p)k}\sum_{n\geq \lfloor t\rfloor}\e^{-n/3}\leq  18\sqrt{(1-p)k}\,\e^{-t/3}.$$
Squaring both sides, the claim follows.
\QED

We finally turn to the construction of an event $E_n=E_{n,p_n,k_n}$ that
satisfies the properties of Lemma~\ref{lem:goodevent}.

{\noindent{\bf Proof of Lemma~\ref{lem:goodevent}:}\hskip10pt} We assume
$p_n\rightarrow 0$ and fix any sequence of integers $k_n\in\N$ satisfying
$k_n\rightarrow\infty$ and $k_np_n\rightarrow 0$ as
$n\rightarrow\infty$. Concerning property $(a)$ in
Lemma~\ref{lem:goodevent}, we first note that since each new-born child is
a mutant with probability $1-p_n$, we have
$$
\P\left(\tau_{k_n}=b_{k_n}^{(p_n)}\right)=(1-p_n)^{k_n-1}=1-o(1).
$$
Moreover, by~\eqref{eq:behaviortau_n}, we also have $\P(\tau_{k_n}<(\ln
k_n)^2)=1-o(1)$, so that the event 
$$
E_n^1=\left\{\tau_{k_n}=b_{k_n}^{(p_n)}\right\}\cap \left\{\tau_{k_n}<(\ln k_n)^2\right\}
$$
has probability $1-o(1)$ as $n$ tends to infinity.

We turn to property $(b)$. Since conditionally on $\tau_{k_n}$, $(Z(t+\tau_{k_n}): t\geq 0)$ is
a Yule process started from $k_n$ individuals, we obtain
from~\eqref{eq:Yule-distancetostart} and Chebycheff's inequality
$$
\P\left(\sup_{t\geq 0}\big|\e^{-t}Z(t+\tau_{k_n})-k_n| >k_n^{2/3}
\right)\leq \frac{4}{k_n^{1/3}}.
$$
In particular, if we let
$$
E_n^2=\left\{(1-k_n^{-1/3})k_n\e^t\,\leq Z(\tau_{k_n}+t)\, \leq
  (1+k_n^{-1/3})k_n\e^t\quad\textup{for all }t\geq 0\right\},
$$
then $\P(E_n^2)\geq 1-4/k_n^{1/3}=1-o(1)$. 

Finally, for property $(c)$ we need control over $(T^{(p)}(t):t\geq 0)$. In
this regard, Lemma~\ref{lem:numberoftypes} in combination with Chebycheff's
inequality shows that conditionally on $\{\tau_{k_n}=b_{k_n}^{(p_n)}\},$
the event
$$
E_n^3=\left\{\left|T^{(p_n)}(\tau_{k_n}+t)
    -\left(k_n+\int_{\tau_{k_n}}^{\tau_{k_n}+t}(1-p_n)Z(r)\dt
      r\right)\right|\leq \e^{5t/6}k_n^{2/3}\quad\textup{for all }t\geq
  0\right\}
$$
has probability at least $1-324/k_n^{1/3}$. On $E_n^2\cap E_n^3$, we note
that for $n$ sufficiently large, recalling that $p_nk_n\rightarrow 0$,
$$
T^{(p_n)}(\tau_{k_n}+t)\geq k_n+(1-p_n)(1-k_n^{-1/3})k_n\int_0^{t}\e^r\dt r -
\e^{5t/6}k_n^{2/3}\geq  (1-3k_n^{-1/3})k_n\e^t, 
$$
and similarly
$$
T^{(p_n)}(\tau_{k_n}+t)\leq (1+3k_n^{-1/3})k_n\e^t.
$$
Setting $$E_n=E_n^1\cap E_n^2\cap E_n^3,$$ we have constructed an event of
probability $1-o(1)$ on which for $n$ sufficiently large, $(a)$, $(b)$ and
$(c)$ in the statement of Lemma~\ref{lem:goodevent} are fulfilled.  \QED
\end{appendix}


\begin{thebibliography}{99}

\bibitem{Aldous} {Aldous, D.}  {\sl Probability approximations via the
    {P}oisson clumping heuristic}. Springer-Verlag, New York, (1989).

\bibitem{AlSp} {Alon, N., Spencer, J.} {\sl The probabilistic method.} Wiley,
  Third Edition (2008).

\bibitem{AtNe} {Athreya, K.B., Ney, P.E.} {\sl Branching Processes.} Dover
  Books on Mathematics (2004). 

\bibitem{Ba} {Baur, E.} {Percolation on random recursive trees.} Preprint
  (2014), available at arXiv:1407.2508. To appear in {\sl Random Struct. Algor.}

\bibitem{BaBe} {Baur, E., Bertoin, J.} {The fragmentation process of an
    infinite recursive tree and Ornstein-Uhlenbeck type processes.}  {\sl
    Electron. J. Probab.}  {\bf 20 (98)} (2015), 1--20.

\bibitem{Be2} {Bertoin, J.} {Sizes of the largest clusters for
    supercritical percolation on random recursive trees.}  {\it Random
    Struct. Algor.} {\bf 44 (1)} (2014), 1098--2418. 

\bibitem{Be4} {Bertoin, J.} {On the non-Gaussian fluctuations of the giant
    cluster for percolation on random recursive trees.}  {\it Electron.
    J. Probab.}  {\bf 19} (2014), no. 24, 1--15.

\bibitem{BeBr} {Bertoin, J., Uribe Bravo, G.} {Supercritical percolation on
  large scale-free random trees.} {\it Ann. Appl. Probab.} {\bf 25 (1)} (2015), 81--103.  

\bibitem{Ber} {Berzunza, G.} {Yule processes with rare mutation and their
  applications to percolation on b-ary trees.}  {\it Electron.  J. Probab.}
  {\bf 20} (2015), 1--23.

\bibitem{Bi} {Billinsgley, P.} {\sl Convergence of Probability Measures.}
  Second edition. Wiley series in Probability and Statistics (1999).

\bibitem{Bo} {Bollob\'as, B.} {\sl Random Graphs.}  Second
  edition. Cambridge University Press (2001).

\bibitem{DoSm} {Dobrow, R. P., Smythe R. T.} {Poisson approximations for
    functionals of random trees.} {\it Random Struct. Algor.}, {\bf
    9 (1-2)} (1996), 79--92.
    
\bibitem{DIMR} {Drmota, M., Iksanov, A., M\"ohle, M. and R\"osler, U}. {A
    limiting distribution for the number of cuts needed to isolate the root
    of a random recursive tree}. {\it Random Struct. Algor.} {\bf
    34 (3)} (2009), 319--336.

\bibitem{GM} {Goldschmidt, C. and Martin, J. B.} {Random recursive trees
    and the {B}olthausen-{S}znitman coalescent.} {\it Electron. J. Probab.}
  {\bf 10} (2005), 718--745.
    
\bibitem{IM} {Iksanov, A. and M{\"o}hle, M.} {A probabilistic proof of a
    weak limit law for the number of cuts needed to isolate the root of a
    random recursive tree.} {\it Electron. Comm. Probab.}  {\bf 12} (2007),
  28--35.

\bibitem{Kl} {Klebaner, F. C.} {\sl Introduction to Stochastic Calculus
    with Applications.} Imperial College Press, 2nd edition (2005).
    
\bibitem{KP} {Kuba, M. and Panholzer, A.} {Multiple isolation of nodes in
    recursive trees.}  {\it Online J. Anal. Comb.} {\bf 9},
  (2014). Available at http://www.math.rochester.edu/ojac/vol9/98.pdf

\bibitem{LC} {Le Cam, L.} {An approximation theorem for the Poisson
    binomial distribution.} {\it Pacific J. Math.} {\bf 10} (1960), 1181--1197. 
    
\bibitem{MM} {Meir, A. and Moon, J. W.} {Cutting down recursive
trees}. {\it Math. Biosci.} {\bf 21} (1974), 173--181.

\bibitem{Mo} {M\"ohle, M.} {The Mittag-Leffler process and a scaling limit
    for the block counting process of the Bolthausen-Sznitman coalescent.}
  {\it Lat. Am. J. Probab. Math. Stat.} {\bf 12 (1)} (2015), 35--53.

\end{thebibliography}
\end{document}